\documentclass[thmsb,final,a4paper,12pt]{article}
\usepackage{amssymb}
\usepackage{amsfonts}
\usepackage{amsmath}
\usepackage{graphicx}
\usepackage{geometry}

\newtheorem{theorem}{Theorem}[section]
\newtheorem{state}{state}[section]

\newtheorem{corollary}[state]{Corollary}

\newtheorem{definition}{Definition}[section]

\newtheorem{lemma}[state]{Lemma}
\newtheorem{notation}{Notation}[section]

\newtheorem{proposition}[state]{Proposition}
\newtheorem{remark}[state]{Remark}

\newenvironment{proof}[1][Proof]{\textbf{#1.} }{\ \rule{0.5em}{0.5em}}

\input{tcilatex}
\begin{document}

\title{Random walk on graphs with regular resistance and volume growth%
\thanks{%
Running head: Random walks on graphs}}
\author{Andr\'{a}s Telcs \\
{\small Department of Computer Science and Information Theory, }\\
{\small Budapest} {\small University of Technology and Economics}\\
{\small telcs@szit.bme.hu}}
\maketitle

\begin{abstract}
In this paper characterizations of graphs satisfying heat kernel estimates
for a wide class of space-time scaling functions are given. The equivalence
of the two-sided heat kernel estimate and the parabolic Harnack inequality
is also shown via the equivalence of the upper (lower) heat kernel estimate
to the parabolic mean value (and super mean value) inequality.
\end{abstract}

\tableofcontents

\section{Introduction}

\setcounter{equation}{0}\label{sintr}The heat propagation through a medium
is determined by its heat capacity and conductance of the media. \ For the
Euclidian space this observation goes back to Einstein. \ The heat diffusion
has been a subject of interest in the discrete and continuous case for
several decades and the fundamental results go back to the classical works
of Aronson \cite{A1}, Davies \cite{D1}, Fabes, Stroock \cite{FS1},\
Grigor'yan, \cite{G1}, Moser \cite{M1},\cite{M2}, Li,Yau \cite{LY1},
Saloff-Coste \cite{SC}, Varadhan \cite{V1}. All these works are confined to
homogeneous spaces. The diffusion in these spaces is typically located
within the distance $\sqrt{t}$, at time $t$, from the starting point. In
other words, the time-space scaling is $\left( time\right) ^{1/2}$ or the
space-time scaling is $\left( distance\right) ^{2}$. The inhomogeneous case
attracts more and more attention of physicist and mathematicians since the
80-s. Geometric and algebraic conditions are relaxed and fractals enriched
the topics. ( For recent results see \cite{CG}, \cite{Gfk}, \cite{Sab}, \cite%
{HS}, \cite{BCK}, \cite{BB} and \cite{BBK}.)

Delmotte has shown in \cite{De} (in the spirit of the results for manifolds
by Saloff-Coste \cite{SC} and Grigor'yan \cite{G1}) for general graphs that
the two-sided Gaussian heat kernel estimate%
\begin{equation}
\frac{c\exp \left( -C\frac{d\left( x,y\right) ^{2}}{n}\right) }{V\left( x,%
\sqrt{n}\right) }\leq \widetilde{p}_{n}\left( x,y\right) \leq \frac{C\exp
\left( -c\frac{d\left( x,y\right) ^{2}}{n}\right) }{V\left( x,\sqrt{n}%
\right) }  \label{GEV2}
\end{equation}%
is equivalent to the parabolic Harnack inequality (with $R=R^{2}$ scaling,
see all the formal definitions below).

In the last two decades several works have been devoted to fractals and
fractal like graphs. \ One of the particular features of these structures is
that the walk (or process) admits the space-time scaling function $R^{\beta
} $ with an exponent $\beta >2.$ \ For the continuous case the equivalence
of the two-sided heat kernel estimate and the parabolic Harnack inequality
with $F\left( x,R\right) =R^{\beta }$ scaling has been shown in {\small %
Hebisch, Saloff-Coste, }\cite{HS} (see also Barlow, Bass \cite{BB} and
Barlow, Bass and Kumagai \ \cite{BBK}). In the graph case, the equivalence
to the two-sided sub-Gaussian estimate (for $\beta >1$) 
\begin{equation}
\frac{c\exp \left[ -C\left( \frac{d^{\beta }\left( x,y\right) }{n}\right) ^{%
\frac{1}{\beta -1}}\right] }{V\left( x,n^{\frac{1}{\beta }}\right) }\leq 
\widetilde{p}_{n}\left( x,y\right) \leq \frac{C\exp \left[ -c\left( \frac{%
d^{\beta }\left( x,y\right) }{n}\right) ^{\frac{1}{\beta -1}}\right] }{%
V\left( x,n^{\frac{1}{\beta }}\right) }  \label{sgebeta}
\end{equation}%
and other conditions was shown in \cite{GT2}. For a wider set of space-time
scaling functions $F\left( x,R\right) =F\left( R\right) $ the corresponding
results were obtained in \cite{t2}.

\ Barlow, Coulhon, Grigor'yan investigated in \cite{BCG} the long time
behavior of the heat kernel on manifolds using volume growth conditions.
Here a more detailed picture will be provided covering on- and off-diagonal
estimates under volume growth and potential theoretic conditions. \ 

Among others Hino, Ram\'{\i}rez \cite{HR1}, Norris \cite{N1} \ and Sturm 
\cite{ST1} (see also references there) studied the heat diffusion in
Dirichlet spaces. \ Their approach uses the intrinsic metric which recovers
the classical Gaussian heat kernel estimate. For us the metric is a priori
given and the space-time scaling might be different from the classical $%
R^{2} $. (For more comment about the difference of the two approaches with
respect of fractals see comments in Section 3.2 of \cite{HR1}.)

The present paper is partly motivated by the works Li and Wang \cite{LW} and
Sung \cite{Su}. We prove in the context of weighted graphs that for a wide
set of scaling functions the heat kernel upper estimate is equivalent to the
parabolic mean value inequality (among others this is shown in \cite{LW} for
the $R^{2}$ scaling). We also show ( inspired by \cite{Su} confined to the $%
R^{2}$ scaling) that some lower estimates are equivalent to the super mean
value inequality. As a consequence, we prove that the conjunction of the
parabolic mean value and super mean value inequality is equivalent to the
two-sided heat kernel estimate and to the parabolic Harnack inequality, as
well.

Recent studies successfully transfer results obtained in continuous setting
to the discrete graph case and vice versa (c.f. \cite{BBK}, \cite{GT2}, \cite%
{GT3} ). For instance in \cite{BBK} the proof of the equivalence of the
parabolic Harnack inequality and two-sided heat kernel estimate (for the $%
R^{\beta }$ scaling) is given for measure metric Dirichlet spaces via to the
graph case where the equivalence is known (c.f. \cite{BB}). We treat the
graph case while we believe that all the arguments and results can be
transferred and are valid for measure metric spaces equipped with a strongly
local, regular symmetric Dirichlet form and with the corresponding diffusion
process.

The aim of the present paper is to relax, as much as possible, the
conditions imposed on the space-time scaling function. We will consider
graphs for which the space-time scaling function $F\left( x,R\right) $ is
not uniform in the center $x$. One may feel that such a generalization is
formal. \ A very simple example shows the opposite (see \cite{t2}), for the
constructed graph neither the volume $V\left( x,R\right) $ \ nor the
space-time scaling function $F\left( x,R\right) $ is uniform in $x\in \Gamma 
$ \ but heat kernel estimates hold.

\ Fractafolds, defined by Strichartz \cite{St}, are the continuous
counterparts of such structures and as it is mentioned above we expect that
the presented results are transferable to continuous spaces and to
fractafolds.

Before we can state our results we need some definitions.

\subsection{Basic definitions}

Let us consider a countable infinite connected graph $\Gamma $. A weight
function $\mu _{x,y}=\mu _{y,x}>0$ is given on the edges $x\sim y.$ This
weight induces a measure $\mu (x)$%
\begin{equation*}
\mu (x)=\sum_{y\sim x}\mu _{x,y},\text{ }\mu (A)=\sum_{y\in A}\mu (y)
\end{equation*}%
on the vertex set $A\subset \Gamma $ and defines a reversible Markov chain $%
X_{n}\in \Gamma $, i.e. a random walk on the weighted graph $(\Gamma ,\mu )$
with transition probabilities 
\begin{align*}
P(x,y)& =\frac{\mu _{x,y}}{\mu (x)}, \\
P_{n}(x,y)& =\mathbb{P}(X_{n}=y|X_{0}=x)\text{ and the corresponding kernel,}
\\
p_{n}(x,y)& =\frac{1}{\mu \left( y\right) }P_{n}(x,y).
\end{align*}%
Let us use the notation$\ \widetilde{p}_{n}=p_{n+1}+p_{n}$. The graph is
equipped with the usual (shortest path length) graph distance $d(x,y)$ and
open metric balls are defined for $x\in \Gamma ,$ $R>0$ as $B(x,R)=\{y\in
\Gamma :d(x,y)<R\}.$ The $\mu -$measure of balls is denoted by 
\begin{equation}
V(x,R)=\mu \left( B\left( x,R\right) \right) .  \label{vdef}
\end{equation}%
For a set $A\subset \Gamma $ the killed random walk is defined by the
transition operator restricted to $c_{0}\left( A\right) \ $(to the set of
functions with support in $A$) and the corresponding transition probability
and kernel is denoted by $P_{n}^{A}\left( x,y\right) $ and $p_{k}^{A}\left(
x,y\right) $.

The (heat) kernel $p_{n}\left( x,y\right) $ is the fundamental solution of
the discrete heat equation on $\left( \Gamma ,\mu \right) :$%
\begin{equation}
\partial _{n}u=\Delta u,  \label{he}
\end{equation}%
where $\partial _{n}u=u_{n+1}-u_{n}$ is the discrete differential operator
with respect of the time and $\Delta =P-I$ is the Laplace operator on $%
\Gamma $.

\begin{definition}
Throughout the paper we will assume that condition $\mathbf{(p}_{0}\mathbf{)}
$ holds, that is, there is a universal $p_{0}>0$ such that for all $x,y\in
\Gamma ,x\sim y$ 
\begin{equation}
\frac{\mu _{x,y}}{\mu (x)}\geq p_{0}.  \label{p0}
\end{equation}
\end{definition}

\begin{definition}
\label{dvd}The weighted graph has the volume doubling\emph{\ }$\left( 
\mathbf{VD}\right) $ property if there is a constant $D_{V}>0$ such that for
all $x\in \Gamma $ and $R>0$%
\begin{equation}
V(x,2R)\leq D_{V}V(x,R).  \label{PD1V}
\end{equation}
\end{definition}

\begin{notation}
For convenience we introduce a short notation for the volume of the annulus: 
$v(x,r,R)=V(x,R)-V(x,r)$ for $R>r>0,x\in \Gamma $.
\end{notation}

\begin{definition}
Now let us consider the exit time%
\begin{equation*}
T_{B(x,R)}=\min \{k:X_{k}\notin B(x,R)\}
\end{equation*}%
from the ball $B(x,R)$ and its mean value 
\begin{equation*}
E_{z}(x,R)=\mathbb{E}(T_{B(x,R)}|X_{0}=z)
\end{equation*}%
and let us use the notation 
\begin{equation*}
E(x,R)=E_{x}(x,R).
\end{equation*}
\end{definition}

\begin{definition}
We will say that the weighted graph $(\Gamma ,\mu )$ satisfies the time
comparison principle $\left( \mathbf{TC}\right) $ if there is a constant $%
C_{T}>1$ such that for all $x\in \Gamma $ and $R>0,y\in B\left( x,R\right) $%
\begin{equation}
\frac{E(x,2R)}{E\left( y,R\right) }\leq C_{T}.  \label{TC}
\end{equation}
\end{definition}

\begin{definition}
We will say that the weighted graph $(\Gamma ,\mu )$ satisfies the weak time
comparison principle $\left( \mathbf{wTC}\right) $ if there is a constant $%
C>1$ such that for all $x\in \Gamma $ and $R>0,y\in B\left( x,R\right) $%
\begin{equation}
\frac{E(x,R)}{E\left( y,R\right) }\leq C.  \label{wTC}
\end{equation}
\end{definition}

\begin{notation}
For a set $A\subset \Gamma $ denote the closure by 
\begin{equation*}
\overline{A}=\left\{ y\in \Gamma :\text{there is an }x\in A\text{ such that }%
x\sim y\right\} ,
\end{equation*}%
the boundary $\partial A=\overline{A}\backslash A$ and $A^{c}=\Gamma
\backslash A.$
\end{notation}

\begin{definition}
A function $h$ is harmonic on a set $A\subset \Gamma $ if it is defined on $%
\overline{A}$ and 
\begin{equation*}
Ph\left( x\right) =\sum_{y}P\left( x,y\right) h\left( y\right) =h\left(
x\right)
\end{equation*}%
for all $x\in A.$
\end{definition}

\begin{definition}
The weighted graph $\left( \Gamma ,\mu \right) $ satisfies the elliptic
Harnack inequality $(\mathbf{H})$ if there is a $C>0$ such that for all $%
x\in \Gamma $ and $R>0$ and for all $u\geq 0$ harmonic functions on $B(x,2R)$
the following inequality holds 
\begin{equation}
\max_{B(x,R)}u\leq C\min_{B(x,R)}u.  \label{H}
\end{equation}
\end{definition}

One can check easily that for any fixed $R_{0}$ for all $R<R_{0}$ the
Harnack inequality follows from $\left( p_{0}\right) $.

\begin{definition}
We define $W_{0}$ to be the set of functions which are candidates to be a
space-time scaling function. In particular$F\in W_{0}$ if $F:\Gamma \times 
\mathbb{N}
\rightarrow 
\mathbb{R}
$ and\newline
1. there are $\beta >1,\beta ^{\prime }>0,c_{F},C_{F}>0$ such that for all $%
R>r>0,x\in \Gamma ,y\in B\left( x,R\right) $%
\begin{equation}
c_{F}\left( \frac{R}{r}\right) ^{\beta ^{\prime }}\leq \frac{F\left(
x,R\right) }{F\left( y,r\right) }\leq C_{F}\left( \frac{R}{r}\right) ^{\beta
},  \label{w1}
\end{equation}%
2. there is a $c>0$ such that for all $x\in \Gamma ,R>0$ 
\begin{equation}
F\left( x,R\right) \geq cR^{2},  \label{w2}
\end{equation}%
3. 
\begin{equation}
F\left( x,R+1\right) \geq F\left( x,R\right) +1  \label{fmonot}
\end{equation}%
for all $R\in 
\mathbb{N}
.$\newline
\newline
Finally $F\in W_{1}$ if $F\in W_{0}$ and $\beta ^{\prime }>1$ holds as well.
\end{definition}

\begin{remark}
\label{rinv}We have by $\left( \ref{fmonot}\right) $ that the function $%
F\left( x,R\right) $ is strictly increasing from $\mathbb{N}$ to $\mathbb{R}$
in the second variable consequently it has the generalized inverse $f:%
\mathbb{N\rightarrow N}:$%
\begin{equation*}
f\left( x,n\right) =\min \left\{ R\in \mathbb{N}:F\left( x,R\right) \geq
n\right\} .
\end{equation*}%
In the whole sequel $f\left( x,n\right) $ is reserved for this inverse.
\end{remark}

The function sets $W_{1}\subset W_{0}$ will play a particular role in the
whole sequel.

Sometimes we will refer to the upper and lower estimate in $\left( \ref{w1}%
\right) $ for $x=y$ as the doubling and the anti-doubling property and in
general, jointly we refer to them as doubling or regularity properties.

\begin{definition}
\label{dpH}We say that $\mathbf{PH}\left( F\right) $, the parabolic Harnack
inequality holds for a function $F$ if \ for the weighted graph $\left(
\Gamma ,\mu \right) $ there is a constant $C>0$ such that for any $x\in
\Gamma ,R,k\geq 0$ and any solution $u\geq 0$ of the heat equation $\left( %
\ref{he}\right) $ on $\mathcal{D}=[k,k+F(x,R)]\times B(x,2R)$ the following
is true. On the smaller cylinders defined by 
\begin{eqnarray*}
\mathcal{D}^{-} &=&[k+\frac{1}{4}F(x,R),k+\frac{1}{2}F(x,R)]\times B(x,R)%
\text{ } \\
\text{and }\mathcal{D}^{+} &=&[k+\frac{3}{4}F(x,R),k+F(x,R))\times B(x,R)
\end{eqnarray*}%
and taking $(n_{-},x_{-})\in \mathcal{D}^{-},(n_{+},x_{+})\in \mathcal{D}%
^{+},$%
\begin{equation}
d(x_{-},x_{+})\leq n_{+}-n_{-}  \label{smdist}
\end{equation}%
the inequality 
\begin{equation*}
u_{n_{-}}(x_{-})\leq C\widetilde{u}_{n_{+}}(x_{+})
\end{equation*}%
holds, where we use the $\widetilde{u}_{n}=u_{n}+u_{n+1}$.
\end{definition}

\begin{definition}
Let us define the "volume" of a space-time cylinder $D=\left[ n,m\right]
\times B\left( x,R\right) $ where $m>n,~R>0$ by 
\begin{equation*}
\nu \left( D\right) =\left[ m-n\right] V\left( x,R\right) .
\end{equation*}
\end{definition}

\begin{definition}
\label{dpmv}We say that $\mathbf{PMV}_{\delta }\left( F\right) $ (the strong
form of) the parabolic mean value inequality$\ $with respect to a function $%
F $ holds on $(\Gamma ,\mu )$ if for fixed constants $0\leq
c_{1}<c_{2}<c_{3}<c_{4}\leq c_{5},0<\delta \leq 1$ there is a $C>1$ such
that for arbitrary $x\in \Gamma $ and $R>0,$ using the notations $F=F\left(
x,R\right) ,B=B\left( x,R\right) ,$ $\mathcal{D}=\left[ 0,c_{5}F\right]
\times B,$ $\mathcal{D}^{-}=\left[ c_{1}F,c_{2}F\right] \times B\left(
x,\delta R\right) ,$ $\mathcal{D}^{+}=\left[ c_{3}F,c_{4}F\right] \times
B\left( x,\delta R\right) $ for any non-negative Dirichlet sub-solution of
the heat equation 
\begin{equation*}
\triangle ^{B}u\geq \partial _{n}u
\end{equation*}%
on $\mathcal{D},$ the inequality 
\begin{equation}
\max_{\mathcal{D}^{+}}u\leq \frac{C}{\nu \left( \mathcal{D}^{-}\right) }%
\sum_{\left( i,y\right) \in \mathcal{D}^{-}}u_{i}(y)\mu (y)\,  \label{PMV}
\end{equation}%
holds.
\end{definition}

\begin{definition}
We will use $PMV\left( F\right) $ if $PMV_{\delta }\left( F\right) $ holds
for $\delta =1.$
\end{definition}

\begin{definition}
We say that (the strong form of) the parabolic super mean value inequality $%
\mathbf{PSMV}\left( F\right) $ holds on $(\Gamma ,\mu )$ with respect to a
function $F$ if there is an $0<\varepsilon <1$ such that for any constants $%
0<c_{1}<c_{2}<c_{3}<c_{4}\leq c_{5},$ $c_{4}-c_{1}<\varepsilon ,$ there are $%
\delta ,$ $c>0$ such that for arbitrary $x\in \Gamma $ and $R>0,$ using the
notations $\ F=F\left( x,R\right) ,B=B\left( x,R\right) ,\mathcal{D}=\left[
0,c_{5}F\right] \times B,$ $\mathcal{D}^{+}=\left[ c_{3}F,c_{4}F\right]
\times B\left( x,\delta R\right) $, $\mathcal{D}^{-}=\left[ c_{1}F,c_{2}F%
\right] \times B\left( x,\delta R\right) $ for any non-negative Dirichlet
super-solution of the heat equation 
\begin{equation*}
\triangle ^{B}u\leq \partial _{n}u
\end{equation*}%
on $\mathcal{D},$ the inequality 
\begin{equation}
\min_{\mathcal{D}^{+}}\widetilde{u}_{k}\geq \frac{c}{\nu \left( \mathcal{D}%
^{-}\right) }\sum_{\left( i,y\right) \in \mathcal{D}^{-}}\widetilde{u}%
_{i}(y)\mu (y)\,  \label{minPH}
\end{equation}%
holds.
\end{definition}

\begin{definition}
We introduce for $A\subset \Gamma $%
\begin{equation*}
G^{A}(y,z)=\sum_{k=0}^{\infty }P_{k}^{A}(y,z),
\end{equation*}%
the local Green function, which is the Green function of the killed walk and
the corresponding Green's kernel as 
\begin{equation*}
g^{A}(y,z)=\frac{1}{\mu \left( z\right) }G^{A}(y,z).
\end{equation*}
\end{definition}

\begin{definition}
The Green kernel may satisfy the following properties. There are $c,C>0$ \
and a function $F$ such that for all $x\in \Gamma ,R>0,A=B\left( x,R\right)
\backslash B\left( x,R/2\right) ,B=B\left( x,2R\right) $%
\begin{equation}
\max_{y\in A}g^{B}\left( x,y\right) \leq C\frac{F\left( x,2R\right) }{%
V\left( x,2R\right) },  \label{gUE}
\end{equation}%
\begin{equation}
\min_{y\in A}g^{B}\left( x,y\right) \geq c\frac{F\left( x,2R\right) }{%
V\left( x,2R\right) }  \label{gLE}
\end{equation}%
If both inequalities hold this fact will be denoted by $g\left( F\right) .$
\end{definition}

\subsection{\protect\bigskip Statement of the results}

The main results of the paper are the following.

\begin{theorem}
\label{tmain}If a weighted graph $(\Gamma ,\mu )$ satisfies $\left(
p_{0}\right) $ and $\left( VD\right) $, then the following statements are
equivalent.

\begin{enumerate}
\item There is an $F\in W_{0}$ such that $g\left( F\right) $ \ is satisfied,

\item $\left( wTC\right) $ and $\left( H\right) $ hold,

\item there is an $F\in W_{0}$ such that the upper estimate $\mathbf{UE}%
\left( F\right) $ holds: there are $C,\beta >1,c>0$ such that for all $%
x,y\in \Gamma ,$ $n>0$%
\begin{equation}
p_{n}\left( x,y\right) \leq \frac{C}{V\left( x,f\left( x,n\right) \right) }%
\exp \left[ -c\left( \frac{F\left( x,d\left( x,y\right) \right) }{n}\right)
^{\frac{1}{\beta -1}}\right]  \label{UE}
\end{equation}%
and furthermore the particular lower estimate $\mathbf{PLE}\left( F\right) $
holds: there are \thinspace $0<c,\delta ,\varepsilon <1$ such that for all
\thinspace $x\in \Gamma ,R>0,B=B\left( x,R\right) ,$ $d\left( x,y\right)
<n\wedge \delta f\left( x,n\right) ,$ $n\leq \varepsilon F\left( x,R\right)
, $%
\begin{equation}
\widetilde{p}_{n}^{B}\left( x,y\right) \geq \frac{c}{V\left( x,f\left(
x,n\right) \right) },  \label{PLE0}
\end{equation}%
where $f\left( x,n\right) $ is the inverse of $F\left( x,R\right) $ in the
second variable,

\item there is an $F\in W_{0}$ such that $PMV\left( F\right) $ and $%
PSMV\left( F\right) $ hold.
\end{enumerate}
\end{theorem}

Further equivalent conditions will be given in Section \ref{swhg}.

\begin{remark}
\label{rLE}The off-diagonal lower estimate $\mathbf{LE}\left( F\right) $
which states that there are $C,\beta ^{\prime }>1,c>0$ such that for all $%
x,y\in \Gamma ,$ $n\geq d\left( x,y\right) $%
\begin{equation}
\widetilde{p}_{n}\left( x,y\right) \geq \frac{c}{V\left( x,f\left(
x,n\right) \right) }\exp \left[ -C\left( \frac{F\left( x,d\left( x,y\right)
\right) }{n}\right) ^{\frac{1}{\beta ^{\prime }-1}}\right]  \label{LE}
\end{equation}%
can be obtained from $\left( VD\right) $ and $PLE\left( F\right) $ if $\beta
^{\prime }>1$ in $\left( \ref{w1}\right) $ using Aronson's classical
chaining argument. This indicates the possibility to obtain two-sided heat
kernel estimate and necessary and sufficient conditions for it.
\end{remark}

\begin{theorem}
\label{tph}If a weighted graph $(\Gamma ,\mu )$ satisfies $\left(
p_{0}\right) $, then the following statements are equivalent:

\begin{enumerate}
\item $\left( VD\right) $ holds \ and there is an $F\in W_{1}$ such that $%
g\left( F\right) $ is satisfied,

\item there is an $F\in W_{1}$ such that the two-sided heat kernel estimate
hold: there are $C,\beta \geq \beta ^{\prime }>1,c>0$ such that for all $%
x,y\in \Gamma ,$ $n\geq d\left( x,y\right) $%
\begin{equation}
c\frac{\exp \left[ -C\left( \frac{F\left( x,d\right) }{n}\right) ^{\frac{1}{%
\beta ^{\prime }-1}}\right] }{V\left( x,f\left( x,n\right) \right) }\leq 
\widetilde{p}_{n}\left( x,y\right) \leq C\frac{\exp \left[ -c\left( \frac{%
F\left( x,d\right) }{n}\right) ^{\frac{1}{\beta -1}}\right] }{V\left(
x,f\left( x,n\right) \right) },  \label{tsE}
\end{equation}%
where we write $d=d\left( x,y\right) $,

\item there is an $F\in W_{1}$ such that $PMV\left( F\right) $ and $%
PSMV\left( F\right) $ hold.

\item there is an $F\in W_{1}$ such that $PH\left( F\right) $ holds.
\end{enumerate}
\end{theorem}

Let us remark that the observation that the pair of the upper and particular
lower estimate is equivalent to the $\beta -$parabolic Harnack inequality
goes back to \cite{HS}.

Our results are presented in the discrete case, and with this limitation
(which we consider not essential) are generalization of several works
devoted to heat kernel estimates and the parabolic Harnack inequality for
scaling function $F\left( x,R\right) =R^{2},R^{\beta }$ or $F\left( R\right) 
$, among others \cite{SC},\cite{G1},\cite{CG},\cite{GT1},\cite{GT2},\cite%
{GT3},\cite{Su}.

The following elements of the present paper are new:

\begin{enumerate}
\item the wide sets $W_{0}$ and $W_{1}$ of space-time scaling functions,

\item condition $g\left( F\right) $ with respect to $F\in W_{i},i=0,1$,

\item the parabolic inequalities with respect to $F\in W_{i},i=0,1$,

\item the proof of the equivalence of the conjunction of the parabolic mean-
and super mean value inequality to the parabolic Harnack inequality,

\item the role of the strong anti-doubling property, $\beta ^{\prime }>1$ is
explained.
\end{enumerate}

The condition $g\left( F\right) $ is the generalization of the corresponding
conditions $\left( G\right) $ in \cite{GT1} and $\left( G_{\beta }\right) $
in \cite{GT2}. \ 

In \cite{Su} partial equivalence (for $F\left( x,R\right) =R^{2}$) was shown
for the parabolic super mean value inequality and the particular heat kernel
lower estimate. Here we prove full equivalence for all $F\in W_{0}$ for a
slightly modified version of the parabolic super mean value inequality. \
This modification allows us not only to show the full equivalence, but also
appropriate for proveing that the conjunction of the parabolic mean and
super \ mean value inequality is equivalent to the parabolic Harnack
inequality. We have not found such a result in the literature even for the
classical case $F\left( x,R\right) =R^{2}$.

In most earlier works the space-homogeneous case $E\left( x,R\right) \simeq
F\left( R\right) $ is considered. \ In these situations $\beta ^{\prime }>1$
follows from the homogeneity \ (and from the other conditions needed for the
heat kernel estimates).

The structure of the paper is the following. Section \ref{sdefs} contains
the basic definitions. Section \ref{sER} recalls the results regarding the
mean exit time. \ Section \ref{swhg} contains the proof of Theorem \ref%
{tmain}, first the Subsection \ref{sue} summarizes a result about the heat
kernel upper estimate, and Subsections \ref{sNDLE} and \ref{saPMV} the lower
estimate. The proof of Theorem \ref{tph} which contains the parabolic
Harnack inequality is given in Section \ref{spPH}. The paper is closed with
a short remark on the homogeneous case.

\section{Definitions and preliminaries}

\setcounter{equation}{0}\label{sdefs}

\subsection{The volume}

\begin{definition}
We will use the inner product with respect to $\mu :$%
\begin{equation*}
\left( f,g\right) _{\mu }=\sum_{x\in \Gamma }f\left( x\right) g\left(
x\right) \mu \left( x\right) .
\end{equation*}
\end{definition}

\begin{remark}
One can show that $\left( VD\right) $ is equivalent to 
\begin{equation*}
\frac{V(x,R)}{V(y,S)}\leq C\left( \frac{R}{S}\right) ^{\alpha },
\end{equation*}%
where $\alpha =\log _{2}D_{V}$ and $d(x,y)\leq R$.
\end{remark}

\begin{remark}
It is easy to show (c.f. \cite{CG}) that the volume doubling property
implies the anti-doubling property: there is an $A_{V}>1$ such that for all $%
x\in \Gamma ,R>0$%
\begin{equation}
2V(x,R)\leq V(x,A_{V}R)  \label{aVD}
\end{equation}%
which is equivalent with the existence of $c,\alpha ^{\prime }>0$ such that
for all $x\in \Gamma ,R>r>0$%
\begin{equation*}
\frac{V\left( x,R\right) }{V\left( x,r\right) }\geq c\left( \frac{R}{r}%
\right) ^{\alpha ^{\prime }}.
\end{equation*}
\end{remark}

\begin{notation}
For two real series $a_{\xi },b_{\xi },\xi \in S$ we shall use the notation $%
a_{\xi }\simeq b_{\xi }$ \ if there is a $C>1$ such that for all $\xi \in S$%
\begin{equation*}
C^{-1}a_{\xi }\leq b_{\xi }\leq Ca_{\xi }.
\end{equation*}
\end{notation}

\begin{remark}
Another direct consequence of $\left( p_{0}\right) $ and $(VD)$ is that 
\begin{equation}
v\left( x,R,2R\right) =V(x,2R)-V(x,R)\simeq V(x,R)  \label{V3}
\end{equation}
\end{remark}

\subsection{Laplacian}

\begin{definition}
The random walk on the weighted graph is a reversible Markov chain and the
Markov operator $P$ \ is naturally defined by%
\begin{equation*}
Pf\left( x\right) =\sum P\left( x,y\right) f\left( y\right) .
\end{equation*}
\end{definition}

\begin{definition}
The Laplace operator on the weighted graph $\left( \Gamma ,\mu \right) $ is
defined simply as 
\begin{equation*}
\Delta =P-I.
\end{equation*}
\end{definition}

\begin{definition}
The Laplace operator with Dirichlet boundary conditions on a finite set $%
A\subset \Gamma $ is defined as%
\begin{equation*}
\Delta ^{A}f\left( x\right) =\left\{ 
\begin{array}{ccc}
\Delta f\left( x\right) & \text{if} & x\in A \\ 
0 & if & x\notin A%
\end{array}%
\right. .
\end{equation*}%
The smallest eigenvalue of $-\Delta ^{A}$ is denoted in general by $\lambda
(A)$ and for $A=B(x,R)$ it is denoted by $\lambda =\lambda (x,R)=\lambda
(B(x,R)).$
\end{definition}

\begin{definition}
The energy or Dirichlet form $\mathcal{E}\left( f,f\right) $ associated to $%
\left( \Gamma ,\mu \right) $ is defined as 
\begin{equation*}
\mathcal{E}\left( f,f\right) =-\left( \Delta f,f\right) _{\mu }=\frac{1}{2}%
\sum_{x,y\in \Gamma }\mu _{x,y}\left( f\left( x\right) -f\left( y\right)
\right) ^{2}.
\end{equation*}
\end{definition}

Using this notation the smallest eigenvalue of $-\Delta ^{A}$ can be defined
by 
\begin{equation}
\lambda \left( A\right) =\inf \left\{ \frac{\mathcal{E}\left( f,f\right) }{%
\left( f,f\right) _{\mu }}:f\in c_{0}\left( A\right) ,f\neq 0\right\}
\label{ldef}
\end{equation}%
as well.

\subsection{The resistance}

\begin{definition}
For any two disjoint sets, $A,B\subset \Gamma ,$ the resistance, $\rho
(A,B), $ is defined as 
\begin{equation*}
\rho (A,B)=\left( \inf \left\{ \mathcal{E}\left( f,f\right)
:f|_{A}=1,f|_{B}=0\right\} \right) ^{-1}
\end{equation*}%
and we introduce 
\begin{equation*}
\rho (x,r,R)=\rho (B(x,r),\Gamma \backslash B(x,R))
\end{equation*}%
for the resistance of the annulus around $x\in \Gamma ,$ with $R>r>0$.
\end{definition}

\begin{definition}
We say that the product of the resistance and volume of the annulus is
uniform in the space if%
\begin{equation}
\rho (x,R,2R)v(x,R,2R)\simeq \rho (y,R,2R)v(y,R,2R).  \label{rrvA}
\end{equation}%
We will refer to this property shortly by $\left( \mathbf{\rho v}\right) $.
\end{definition}

\begin{lemma}
\label{crav>l2}For all weighted graphs, $x\in \Gamma ,R>r>0$%
\begin{equation}
\rho (x,r,R)v(x,r,R)\geq (R-r)^{2},  \label{rv>2}
\end{equation}
\end{lemma}

\begin{proof}
For the proof see \cite{ter}.
\end{proof}

\begin{definition}
The resistance lower estimate $\mathbf{RLE}\left( F\right) $ holds for a
function $F$ if there is a $c>0$ such that for all $x\in \Gamma ,R>0$%
\begin{equation}
\rho \left( x,R,2R\right) \geq c\frac{F\left( x,2R\right) }{V\left(
x,2R\right) }.  \label{RLE}
\end{equation}
\end{definition}

\begin{definition}
The anti-doubling property $\left( \mathbf{aD\rho v}\right) $ is satisfied
for $\rho v$ if there are $c,\beta ^{\prime }>0$ such that for all $x\in
\Gamma ,R>r>0$%
\begin{equation}
\frac{\rho \left( x,R,2R\right) v\left( x,R,2R\right) }{\rho \left(
x,r,2r\right) v\left( x,r,2r\right) }\geq c\left( \frac{R}{r}\right) ^{\beta
^{\prime }}.  \label{adrrv}
\end{equation}
\end{definition}

\subsection{The mean exit time}

Let us introduce the exit time $T_{A}$ from a set $A\subset \Gamma $.

\begin{definition}
The exit time from a set $A$ is defined as 
\begin{equation*}
T_{A}=\min \{k:X_{k}\in \Gamma \backslash A\},
\end{equation*}%
its expected value is denoted by 
\begin{equation*}
E_{x}(A)=\mathbb{E}(T_{A}|X_{0}=x)
\end{equation*}%
and furthermore leus us write%
\begin{equation*}
\overline{E}\left( x,R\right) =\max_{y\in B\left( x,R\right) }E_{y}\left(
B\left( x,R\right) \right)
\end{equation*}
\end{definition}

In this section we introduce some properties of the mean exit time which
will play crucial role in the whole sequel. First of all it is immediate that%
\begin{equation*}
E\left( x,1\right) \geq 1
\end{equation*}%
and we for $R\in \mathbb{N}$%
\begin{equation}
E\left( x,R+1\right) \geq E\left( x,R\right) +1.  \label{monot}
\end{equation}

\begin{remark}
We have by $\left( \ref{monot}\right) $ that the function $E\left(
x,R\right) $ is strictly increasing from $\mathbb{N}$ to $\mathbb{R}$ in the
second variable consequently it has the generalized inverse $e:\mathbb{%
N\rightarrow N}$: 
\begin{equation*}
e\left( x,n\right) =\min \left\{ R\in \mathbb{N}:E\left( x,R\right) \geq
n\right\} .
\end{equation*}
\end{remark}

\begin{remark}
It is easy to see that $\left( TC\right) $ is equivalent to the existence of
constants $C,\beta \geq 1$ for which 
\begin{equation}
\frac{E(x,R)}{E(y,S)}\leq C\left( \frac{R}{S}\right) ^{\beta },
\label{tcbeta}
\end{equation}%
for all $y\in B\left( x,R\right) ,R\geq S>0$.
\end{remark}

\begin{definition}
The local sub-Gaussian upper exponent, with respect to a function $F\left(
x,R\right) $ is $k=k_{x}(n,R)\geq 1,$ it is defined as the maximal integer
for which 
\begin{equation}
\frac{n}{k}\leq q\min\limits_{y\in B(x,R)}F(y,\left\lfloor \frac{R}{k}%
\right\rfloor )  \label{iter}
\end{equation}%
or $k=1$ by definition if there is no appropriate $k$. Here $q>0$ is a small
fixed constant (c.f. \cite{TD},$q<\min \left\{ 1/16,c_{F}p_{0}/C_{F}\right\} 
$).
\end{definition}

\begin{definition}
Let $n\geq l_{x}=l_{x}(n,R)\geq 1$ be the minimal integer for which 
\begin{equation}
\frac{n}{l}\geq CF(x,\left\lceil \frac{R}{l}\right\rceil )
\end{equation}%
or $l=n$ by definition if there is no appropriate $l$. The constant $C$ \
will be specified later.
\end{definition}

\begin{definition}
The local sub-Gaussian lower exponent $l\left( n,R,A\right) $ with respect
to a function $F\left( x,R\right) $ for $A\subset \Gamma $ is the maximal
integer $l$ for which%
\begin{equation}
\frac{n}{l}\geq C\max_{z\in A}F(z,\left\lceil \frac{R}{l}\right\rceil ).
\end{equation}
\end{definition}

\begin{definition}
The global sub-Gaussian exponent $m=m\left( n,R\right) $ is defined as the
maximal integer for which 
\begin{equation}
\frac{n}{m}\leq q\min_{y\in \Gamma }E(y,\left\lfloor \frac{R}{m}%
\right\rfloor )
\end{equation}%
or $m=1$ by definition if there is no appropriate $m$.
\end{definition}

\begin{definition}
The mean exit time is uniform in the space if there is a function $F$ such
that%
\begin{equation}
E\left( x,R\right) \simeq F\left( R\right) .  \label{E}
\end{equation}%
This property will be referred to by $\left( E\right) $.
\end{definition}

The definition $m\left( n,R\right) $ is prepared for the particular case
when $E\left( y,R\right) \simeq E\left( x,R\right) \simeq F\left( R\right) $
i.e. $E$ is basically independent of $x,y\in \Gamma .$

\begin{definition}
\label{dv1}We define $V_{1}$ to be the a set of functions such that $F\in
V_{0}$ if $F:\Gamma \times 
\mathbb{N}
\rightarrow 
\mathbb{R}
$ and there are $\beta ^{\prime }>1,c_{F}>0$ such that for all $R>r>0,x\in
\Gamma ,y\in B\left( x,R\right) $%
\begin{equation}
c_{F}\left( \frac{R}{r}\right) ^{\beta ^{\prime }}\leq \frac{F\left(
x,R\right) }{F\left( y,r\right) }.  \label{bb1}
\end{equation}
\end{definition}

\bigskip

\ 

\subsection{Mean value inequalities}

\label{smv}

\begin{definition}
The elliptic mean value inequality $\left( \mathbf{MV}\right) $ holds if
there is a $C>0$ such that for all $x\in \Gamma ,R>0$ and for all $u\geq 0$
harmonic functions on $B=B\left( x,R\right) $%
\begin{equation}
u\left( x\right) \leq \frac{C}{V\left( x,R\right) }\sum_{y\in B}u\left(
y\right) \mu \left( y\right) ,  \label{MV}
\end{equation}
\end{definition}

\begin{remark}
Let us recognize that in the definition of $PMV$ $\delta $ is a "free"
parameter, while in the definition of $PSMV\left( F\right) $ it depends on $%
\varepsilon $ and $c_{i},i=1..5$. Let us also observe that $c_{i}$-s are
subject of the restriction $c_{4}-c_{1}<\varepsilon $ .
\end{remark}

\begin{remark}
One should note that the definition of the parabolic mean value inequality
is slightly different from that is given in \cite{t2}. \ There it is stated
for Dirichlet solutions, here we have it for arbitrary Dirichlet
sub-solutions. It is easy to see that the extended definitions fits into
Theorem \ref{tLDUE+}. On one hand solutions are sub-solutions, on the other
hand, the proof of the implication $UE\Longrightarrow PMV$ of Theorem \ref%
{tLDUE+} follows word by word for sub-solutions.
\end{remark}

\begin{remark}
The condition $\left( \ref{smdist}\right) $ in the definition of the
parabolic Harnack inequality is needed in order to have a path (with nonzero
probability ) of length no more than $n_{+}-n_{-}$ between $x_{-}$ and $%
x_{+} $.\ One can eliminate this restriction if the parabolic Harnack
inequality is considered only for large enough $R$-s. \ The condition $%
n_{+}-n_{-}\geq d(x_{-},x_{+})$ is satisfied if $\left( c_{3}-c_{2}\right)
F\left( x,R\right) \geq cR^{2}>4R$ which holds if $R>R_{0}$. \ Such an $%
R_{0} $ depends only on the constants (c.f. \cite{ter}). In order to avoid
lengthy technical discussion we may assume $R>R_{0}$ in all these
situations. \ If other is not stated, the corresponding inequalities for $%
R\leq R_{0}$ follow from $\left( p_{0}\right) $.
\end{remark}

\section{Properties of the mean exit time}

\setcounter{equation}{0}\label{sER}

In this section we recall some results from \cite{ter} which describe the
behavior of the mean exit time. If this is not mentioned otherwise, the
statements and proofs can be found in \cite{ter}. The first one is the
Einstein relation.

\begin{theorem}
\label{tLER}If $\left( p_{0}\right) ,(VD),\left( H\right) $ \ and one of the
conditions $\left( wTC\right) ,\left( aD\rho v\right) ,RLE\left( E\right) $
or $\rho v\in W_{0}$ hold, then $\left( ER\right) $, the Einstein relation 
\begin{equation}
E(x,2R)\simeq \rho (x,R,2R)v(x,R,2R)  \label{ER}
\end{equation}%
holds, furthermore%
\begin{equation*}
\lambda ^{-1}\left( x,R\right) \simeq E\left( x,R\right) \simeq \overline{E}%
\left( x,R\right)
\end{equation*}%
and 
\begin{equation*}
E\left( x,R\right) \in W_{0}.
\end{equation*}
\end{theorem}

\begin{theorem}
\label{tER1}If $\left( p_{0}\right) ,(VD),\left( TC\right) $ hold, then the
Einstein relation 
\begin{equation}
E(x,2R)\simeq \rho (x,R,2R)v(x,R,2R)
\end{equation}%
holds, furthermore%
\begin{equation*}
\lambda ^{-1}\left( x,R\right) \simeq \overline{E}\left( x,R\right) \simeq
E\left( x,R\right) \in W_{0}.
\end{equation*}
\end{theorem}

The properties of the inverse function $e$ and properties of $E$ are linked
as the following evident lemma states.

\begin{lemma}
The following statements are equivalent \newline
\newline
1. There are $C,c>0,\beta \geq \beta ^{\prime }>0$ such that for all $x\in
\Gamma ,R\geq r>0,$ $y\in B\left( x,R\right) $%
\begin{equation}
c\left( \frac{R}{r}\right) ^{\beta ^{\prime }}\leq \frac{E\left( x,R\right) 
}{E\left( y,r\right) }\leq C\left( \frac{R}{r}\right) ^{\beta };  \label{Eb}
\end{equation}%
2. There are $C,c>0,\beta \geq \beta ^{\prime }>0$ such that for all $x\in
\Gamma ,n\geq m>0,$ $y\in B\left( x,e\left( x,n\right) \right) $%
\begin{equation}
c\left( \frac{n}{m}\right) ^{1/\beta }\leq \frac{e\left( x,n\right) }{%
e\left( y,m\right) }\leq C\left( \frac{n}{m}\right) ^{1/\beta ^{\prime }}.
\label{eb}
\end{equation}
\end{lemma}

\bigskip

\begin{lemma}
\label{k> copy(1)} If $F\in W_{0}$, then for$\ k_{x}(n,R)$ defined in $%
\left( \ref{iter}\right) $ 
\begin{equation}
k_{x}(n,R)+1\geq c\left( \frac{F(x,R)}{n}\right) ^{\frac{1}{\beta -1}}\text{%
\ }  \label{k>e/n}
\end{equation}%
for all $x\in \Gamma ,R,n>0$ for fixed $c>0,\beta >1.$ In addition, if $F\in
W_{1}$ is assumed, then%
\begin{equation}
l_{x}\left( n,R\right) -1\leq C\left( \frac{F(x,R)}{n}\right) ^{\frac{1}{%
\beta ^{\prime }-1}}\text{\ }.  \label{l<e/n}
\end{equation}
\end{lemma}

\begin{proof}
The statement follows from the regularity properties of $F$ easily, $\beta
>1 $ is ensured by $\beta \geq 2$ (see $\left( \ref{w2}\right) $) and $\beta
^{\prime }>1$ by the assumption.
\end{proof}

\begin{lemma}
\label{klm}If $\left( E\right) $ is satisfied\ and \ $E\in W_{0}$ then%
\begin{equation*}
k_{x}\left( n,R\right) \simeq l_{x}\left( n,R\right) \simeq m\left(
n,R\right) .
\end{equation*}
\end{lemma}

Since this fact is not used in the proof of the main results, the elementary
proof is omitted (for some hints see \cite{TD}).

\begin{remark}
Let us mention here that under $\left( p_{0}\right) ,\left( VD\right) $ and $%
\left( H\right) $ the uniformity of the mean exit time in the space: 
\begin{equation*}
E\left( x,R\right) \simeq F\left( R\right)
\end{equation*}%
ensures that $E$ satisfies the left hand side of $\left( \ref{w1}\right) $
with a $\beta ^{\prime }>1$ (c.f. \cite{ter}). This explains that in the
"classical" cases, when $\left( E\right) $ holds one should not assume $%
\beta ^{\prime }>1$, it follows from the conditions (see also \cite{t2}
Theorem 4.13).
\end{remark}

\begin{lemma}
\label{lminE<rv}For $\left( \Gamma ,\mu \right) $ for all $x\in \Gamma ,R>0$%
\begin{equation}
\min_{z\in \partial B\left( x,\frac{3}{2}R\right) }E\left( z,R/2\right) \leq
\rho \left( x,R,2R\right) v\left( x,R,2R\right) .  \label{E<rrv}
\end{equation}
\end{lemma}

\begin{corollary}
\label{cg&H}Under $\left( p_{0}\right) $ and $\left( VD\right) $%
\begin{equation*}
\exists F\in W_{0}:g\left( F\right) \Leftrightarrow \left( H\right) +\left(
ER\right) .
\end{equation*}
\end{corollary}

\begin{proof}
The implication $\Longleftarrow $ was shown in  \cite{ter}. We also know
that\  $\left( p_{0}\right) ,\left( VD\right) ,\left( H\right) $ and $\left(
ER\right) $ implies $\left( TC\right) $ hence by Theorem \ref{tER1} $E\in
W_{0}$. The reverse implication needs some additional arguments. We know
again from \cite{ter} that $g\left( F\right) $ implies $\left( H\right) $.
We show here the implication $g\left( F\right) \Rightarrow \left( ER\right) $
under $\left( p_{0}\right) ,\left( VD\right) $ and $\left( H\right) $. \
That  needs some care. Let us assume that $r_{i}=2^{i},r_{n-1}<2R\leq
r_{n},B_{i}=B\left( x,r_{i}\right) ,,A_{i}=B_{i}\backslash
B_{i-1},V_{i}=V\left( x,r_{i}\right) .$ \ In \cite{GT2} Section 4.3 it is
derived using $\left( p_{0}\right) ,\left( VD\right) $ and $\left( H\right) $%
\ that 
\begin{equation*}
E\left( x,2R\right) \leq C\sum_{i=0}^{n-1}V_{i+1}\rho \left(
x,r_{i},r_{i+1}\right) .
\end{equation*}%
Now we use a consequence of $\left( H\right) $: 
\begin{equation*}
\rho \left( x,r_{i},r_{i+1}\right) \leq C\max_{y\in
A_{i+1}}g^{B_{i+1}}\left( x,y\right) 
\end{equation*}%
(see for instance \cite{ter} Section 4.or \cite{B1}) to obtain%
\begin{eqnarray*}
E\left( x,2R\right)  &\leq &C\sum_{i=0}^{n-1}V_{i+1}\max_{y\in
A_{i+1}}g^{B_{i+1}}\left( x,y\right)  \\
&\leq &C\sum_{i=0}^{n-1}F\left( x,r_{i+1}\right) \leq CF\left(
x,r_{n}\right) \sum_{i=0}^{n-1}2^{-i\beta ^{\prime }} \\
&\leq &CF\left( x,2R\right) ,
\end{eqnarray*}%
where $\left( \ref{gUE}\right) $ was used to get the second inequality.

On the other hand, from $\left( \ref{gLE}\right) $ one obtains 
\begin{eqnarray*}
c\frac{F\left( x,2R\right) }{V\left( x,2R\right) }\left( V\left( x,R\right)
-V\left( x,R/2\right) \right)  &\leq &\min_{y\in B\left( x,R\right)
\backslash B\left( x,R/2\right) }g^{B}\left( x,y\right) \sum_{z\in B\left(
x,R\right) \backslash B\left( x,R/2\right) }\mu \left( z\right)  \\
&\leq &\sum_{z\in B\left( x,R\right) \backslash B\left( x,R/2\right)
}g^{B}\left( x,z\right) \mu \left( z\right) \leq E\left( x,2R\right) ,
\end{eqnarray*}%
this means that 
\begin{equation*}
cF\left( x,2R\right) \leq E\left( x,2R\right) 
\end{equation*}%
consequently, $F\simeq E,$ $E\in W_{0}$ and $\left( TC\right) $ is
satisfied. \ Finally by Theorem \ref{tER1} the conditions $\left(
p_{0}\right) ,\left( VD\right) $ and $\left( TC\right) $ imply $\left(
ER\right) $.
\end{proof}

\begin{corollary}
\label{c*}Assume that $\left( \Gamma ,\mu \right) $ satisfies $\left(
p_{0}\right) ,\left( VD\right) $ and $\left( H\right) $, then%
\begin{equation}
\left( wTC\right) \Leftrightarrow \left( aD\rho v\right) \Leftrightarrow
\left( TC\right) \Leftrightarrow \left( ER\right) \Leftrightarrow RLE\left(
E\right) \overset{\left( i\right) }{\Leftrightarrow }g\left( F\right) ,
\label{*}
\end{equation}%
where $F\in W_{0}$ is a consequence in the direction $\overset{\left(
i\right) }{\Longrightarrow }$ and assumption for $\overset{\left( i\right) }{%
\Longleftarrow }$.
\end{corollary}

\begin{proof}
Exept the last implications the statement was shown in \cite{ter} while the
last one is just Corollary \ref{cg&H}.
\end{proof}

\begin{remark}
\label{rretutn}Let us remark here that as a side result it follows that $%
RLE\left( E\right) $ or $g\left( F\right) $ for $F\in W_{0}$ implies $\rho
v\simeq F$ and $E\simeq F$ as well.
\end{remark}

\section{Temporal regularity and heat kernel estimates}

\setcounter{equation}{0}\label{swhg}In this section we prove Theorem \ref%
{tmain} in an extended form. We have seen in Corollary \ref{c*} that under
the conditions $\left( p_{0}\right) ,\left( VD\right) $ and $\left( H\right) 
$ 
\begin{equation}
\left( wTC\right) \Leftrightarrow \left( aD\rho v\right) \Leftrightarrow
\left( TC\right) \Leftrightarrow \left( ER\right) \Leftrightarrow RLE\left(
E\right) .  \label{**}
\end{equation}%
Let $\left( \ast \right) $ denote any of the equivalent conditions. \ Using
this convention we can state the extension of Theorem \ref{tmain} as follows.

\begin{theorem}
\label{tmainplussz}If a weighted graph $(\Gamma ,\mu )$ satisfies $\left(
p_{0}\right) $ and $\left( VD\right) $, then the following statements are
equivalent:

\begin{enumerate}
\item there is an $F\in W_{0}$ such that $g\left( F\right) $ is satisfied,

\item $\left( H\right) $ and $\left( \ast \right) $ hold,

\item there is an $F\in W_{0}$ such that $UE\left( F\right) $ and $PLE\left(
F\right) $ are satisfied,

\item there is an $F\in W_{0}$ such that $PMV\left( F\right) $ and $%
PSMV\left( F\right) $ are satisfied.
\end{enumerate}
\end{theorem}

The proof \ of Theorem \ref{tmainplussz} contains two autonomous results.
The first one states that the upper estimate is equivalent to the parabolic
mean value inequality, the second one states that the particular lower
estimate is equivalent to the parabolic super mean value inequality. The
return route from $5$ to $1$ and $2$ is based on the Einstein relation, on a
potential theoretic results from $\cite{ter}$ and a modification of the
return route developed in \cite{TD}. The proof of $2\Longrightarrow 3$
generalizes methods of \cite{GT1} and \cite{GT2}.

Let us emphasize the importance of the condition $\left( aD\rho v\right) $
in $\left( \ref{**}\right) $. It is a condition on the volume and
resistance, no assumption of stochastic nature is involved so the result is
in the spirit of Einstein's observation on the heat propagation. \ These
conditions in conjunction with $\left( VD\right) $ \ and $\left( H\right) $
provide the characterization of the heat kernel estimates in terms of volume
and resistance properties. \ Of course the elliptic Harnack inequality is
not easy to verify. Meanwhile we learn from $g\left( F\right) $ that the
main properties ensured by the elliptic Harnack inequality are that the
equipotential surfaces of the local green kernel $g^{B\left( x,R\right) }$
are basically spherical and the potential growth is regular (c.f. \cite{ter}%
).

\subsection{The upper estimate}

\label{sue}This section provides the upper bound part of the implication $%
2\Longrightarrow 3$ of Theorem \ref{tmain}. In details%
\begin{equation}
\left. 
\begin{array}{c}
\left( VD\right) \\ 
\left( TC\right) \\ 
\left( H\right)%
\end{array}%
\right\} \Longrightarrow DUE\left( E\right)
\end{equation}%
and under $\left( VD\right) $ and $\left( TC\right) $

\begin{equation*}
DUE\left( E\right) \Longleftrightarrow UE\left( E\right) \Longleftrightarrow
PMV\left( E\right)
\end{equation*}%
In particular the parabolic mean value inequality is shown to be equivalent
to the upper estimate and the other conditions. This result has been proved
in \cite{t2}:

\begin{theorem}
\label{tLDUE+}For a weighted graph $(\Gamma ,\mu )$ if $\left( p_{0}\right)
,(VD),\left( TC\right) $ conditions hold, then the following statements are
equivalent:\ 

\begin{enumerate}
\item The local diagonal upper estimate $\mathbf{DUE}\left( E\right) $
holds; there is a $C>0$ such that for all $x\in \Gamma ,$ $n>0$ 
\begin{equation}
p_{n}(x,x)\leq \frac{C}{V(x,e(x,n))},  \label{LDUE}
\end{equation}

\item the upper estimate $UE\left( E\right) $ holds: there are $C,\beta
>1,c>0$ such that for all $x,y\in \Gamma ,$ $n>0$%
\begin{equation*}
p_{n}\left( x,y\right) \leq \frac{C}{V\left( x,e\left( x,n\right) \right) }%
\exp \left[ -c\left( \frac{E\left( x,d\left( x,y\right) \right) }{n}\right)
^{\frac{1}{\beta -1}}\right] ,
\end{equation*}

\item the parabolic mean value inequality, $PMV\left( E\right) $ holds,

\item the mean value inequality, $\left( MV\right) $ holds,
\end{enumerate}
\end{theorem}

\begin{corollary}
\label{cDUE}If $\left( \Gamma ,\mu \right) $ satisfies $\left( p_{0}\right) $%
, then%
\begin{equation*}
\left. 
\begin{array}{c}
\left( VD\right) \\ 
\left( wTC\right) \\ 
\left( H\right)%
\end{array}%
\right\} \Longrightarrow UE\left( E\right) .
\end{equation*}
\end{corollary}

\begin{proof}
We know from Theorem \ref{tER1} that $E\in W_{0}$. The statement follows
from Corollary \ref{c*} and Theorem \ref{tLDUE+} since the elliptic Harnack
inequality, $\left( H\right) $ implies the\ elliptic mean value inequality $%
\left( MV\right) $.
\end{proof}

\begin{remark}
The equivalences in $\left( \ref{diagrDUE}\right) $ for $F\in W_{0}$ instead
of $E$ follow from the same proofs given in \cite{t2} for $E$ (see also
Corollary \ref{c*} and Remark \ref{rretutn})
\end{remark}

\begin{remark}
Conequently the implication of the upper bound part in Theorem \ref%
{tmainplussz} (and Theorem \ref{tmain}) $2.\Rightarrow 3.$ is shown.
\end{remark}

\begin{remark}
\label{rpr23}The implication for an $F\in W_{0}$%
\begin{equation*}
DUE\left( F\right) \Longleftrightarrow UE\left( F\right) \Rightarrow
PMV\left( F\right) 
\end{equation*}%
in particular $DUE\left( F\right) \Rightarrow PMV\left( F\right) $ can be
shown repeating step by step the proof given for the particular space-time
scaling function $E$. \ The full proof is spelled out in \cite[Theorem 8.6]%
{ln}. \ This gives the proof of the upper bound part of the implication of
Theorem \ref{tmainplussz} $3.\Rightarrow 4$. 
\end{remark}

\subsection{The near diagonal lower estimate}

\label{sNDLE}In this section we give a lower estimate for the Dirichlet heat
kernel and for the global one.

\begin{definition}
The near diagonal lower estimate, $NDLE\left( F\right) $ holds with respect
to a function $F$ if there are $c,\delta >0$ such that for all \thinspace $%
x,y\in \Gamma ,n>0,d\left( x,y\right) <\delta f\left( x,n\right) \wedge n$%
\begin{equation}
\widetilde{p}_{n}\left( x,y\right) \geq \frac{c}{V\left( x,f\left(
x,n\right) \right) },  \label{NDLE}
\end{equation}%
where $f$ is the existing inverse of $F$ in the second variable, defined in
Remark \ref{rinv}.
\end{definition}

\begin{remark}
It is clear that $PLE\left( F\right) $ \ implies $NDLE\left( F\right) $. It
is also known (c.f.\cite{Gfk}) that $NDLE\left( F\right) $ and $UE\left(
F\right) \ $implies $PLE\left( F\right) $ if $F=R^{2}$, the same proof works 
$F\in W_{0}$.\ 
\end{remark}

\begin{theorem}
\label{tPLE} For weighted graphs 
\begin{equation*}
(p_{0})+\left( VD\right) +\left( TC\right) +(H)\Longrightarrow PLE\left(
E\right)
\end{equation*}
\end{theorem}

The proof closely follows the steps of the corresponding proof given for the
case $E\left( x,R\right) \simeq R^{\beta }$ in \cite{GT1} therefore it is
omitted.

\begin{remark}
From the regularity of $V$\ and $F$\ (and $f$) it is immediate that $%
PLE\left( F\right) $ is equivalent with the slightly stronger form with $%
\delta ^{\prime }=\delta /2:$ for all $x\in \Gamma ,R>0,B=B\left( x,R\right)
,n\leq \varepsilon F\left( x,R\right) ,$ $y,z\in B\left( x,\delta ^{\prime
}f\left( x,n\right) \right) $%
\begin{equation}
\widetilde{p}_{n}^{B}\left( y,z\right) \geq \frac{c}{V\left( y,f\left(
y,n\right) \right) }.  \label{sPLE}
\end{equation}
\end{remark}

\subsection{The parabolic super mean value inequality}

\setcounter{equation}{0}\label{saPMV}In this section the equivalence of the
particular lower estimate and a kind of converse of the parabolic mean value
inequality is shown. \ The partial equivalence for the classical ($F\left(
x,R\right) =$ $R^{2}$) and continuous situation was shown in \cite{Su}. \
Here the generalization to the present settings is provided.

For technical reasons we use some specific constants, like $\varepsilon
,\delta $ from $PLE,$ $c_{F},C_{F}$ from the definition of the set of
scaling functions $W_{0}$ (in $\left( \ref{w1}\right) $).

In this section we show that for an $F\in W_{0}$%
\begin{equation*}
PLE\left( F\right) \Longleftrightarrow PSMV\left( F\right)
\end{equation*}%
that is, the following theorem holds.

\begin{theorem}
\label{tLNLE} For the weighted graph $\left( \Gamma ,\mu \right) $ assume $%
\left( p_{0}\right) $ and $\left( VD\right) $. Then for an $F\in W_{0}$ \ $%
PLE\left( F\right) $ holds if and only if $PSMV\left( F\right) $ holds as
well.
\end{theorem}

\begin{proof}
The proof follows the main steps of \cite{Su}. We know that $PLE\left(
F\right) $ implies \ the slightly stronger version $\left( \ref{sPLE}\right) 
$ that is there are $\delta ,\varepsilon >0$ and $c>0$%
\begin{equation}
p_{m}^{B}(y,z)+p_{m+1}^{B}(y,z)\geq \frac{c}{V(y,f(y,m))}
\end{equation}%
holds, provided that $y,z\in B\left( x,r\right) ,$ where $r=\delta f(x,m)/2,$
and $m\leq \varepsilon F=\varepsilon F(x,R)$. For the super-solution $u$ we
have that for all $c_{3}F\leq k\leq c_{4}F,c_{1}F\leq i\leq c_{2}F$%
\begin{equation*}
\widetilde{u}_{k}\left( y\right) \geq \sum_{z\in B\left( x,R\right) }%
\widetilde{p}_{k-i}^{B}(y,z)\mu \left( z\right) u_{i}\left( z\right) .
\end{equation*}%
In order to use $\left( \ref{sPLE}\right) $, we choose 
\begin{equation}
\delta ^{\ast }=\frac{1}{C_{F}}\left( c_{3}-c_{2}\right) ^{\frac{1}{\beta
^{\prime }}}\delta /2  \label{d*}
\end{equation}%
which ensures that $y,z\in B\left( x,r\right) $ if $r=\delta ^{\ast }R$.
From the condition $c_{4}-c_{1}\leq \varepsilon $ it follows that $k-i\leq
\varepsilon F\left( x,R\right) $ is satisfied and $PLE\left( F\right) $ can
be applied:%
\begin{eqnarray*}
\widetilde{u}_{k}\left( y\right) &\geq &\sum_{y\in B\left( x,\delta ^{\ast
}R\right) }\widetilde{p}_{k-i}^{B}(y,z)\mu \left( z\right) u_{i}\left(
z\right) \\
&\geq &\frac{c}{V(x,f\left( x,k-i\right) )}\sum_{y\in B\left( x,\delta
^{\ast }R\right) }\mu \left( y\right) u_{i}\left( y\right) .
\end{eqnarray*}%
Now let us sum for $c_{1}F\leq i\leq c_{2}F$ and divide by $\left(
c_{2}-c_{1}\right) F$ \ to obtain 
\begin{eqnarray*}
\widetilde{u}_{k}\left( y\right) &\geq &\frac{c}{F\left( x,R\right) }%
\sum_{i=c_{1}F}^{c_{2}F}\frac{1}{V(x,f\left( x,k\right) )}\sum_{y\in B\left(
x,\delta ^{\ast }R\right) }\mu \left( y\right) u_{i}\left( z\right) \\
&\geq &\frac{c}{V\left( x,R\right) F\left( x,R\right) }%
\sum_{i=c_{1}F}^{c_{2}F}\sum_{y\in B\left( x,\delta ^{\ast }R\right) }\mu
\left( y\right) u_{i}\left( z\right) .
\end{eqnarray*}

Now we prove the reverse implication $PSMV\left( F\right) \Longrightarrow
PLE\left( F\right) $ by applying $PSMV$ twice.

\begin{enumerate}
\item Denote $\varepsilon ,\delta ,c_{i}$ the constants in $PSMV\left(
F\right) $, $c_{i}$ will be specified later (which determines $\delta $ as
well), furthermore $F=F\left( p,R\right) ,$ $B=B\left( p,R\right) ,$ $%
r_{1}=R/8,F_{1}=F\left( p,r_{1}\right) ,D_{1}=B\left( p,\delta r_{1}\right)
, $ $m=c^{\prime }F_{1},c_{1}\leq c^{\prime }\leq c_{2}$ and $D=B\left(
p,\delta \frac{R}{4}\right) $. $\,$Let us define 
\begin{equation*}
u_{n}\left( y\right) =\left\{ 
\begin{array}{ccc}
\sum_{z\in D}\widetilde{p}_{n-m}^{B}\left( y,z\right) \mu \left( z\right) & 
\text{if} & n>m \\ 
&  &  \\ 
1 & \text{if} & n\leq m%
\end{array}%
\right. .
\end{equation*}%
This is a solution on $D\times \left[ 0,\infty \right] $ of 
\begin{equation*}
P^{B}u_{n}=u_{n+1}
\end{equation*}%
and $u_{n}\geq 0.$ \ \ From the $PSMV\left( F\right) $ it follows that 
\begin{equation*}
u_{k}\left( x\right) \geq \frac{c}{V\left( x,\delta r_{1}\right) F_{1}}%
\sum_{i=c_{1}F_{1}}^{c_{2}F_{1}}\sum_{w\in D_{1}}\mu \left( w\right) 
\widetilde{u}_{i}\left( w\right)
\end{equation*}%
provided that, $x\in D_{1}$ and $c_{3}F_{1}<k<c_{4}F_{1}.$ From the
definition of $u_{n},$ $\left( VD\right) $ and $F\in W_{0}$ it follows that
\ 
\begin{equation*}
u_{k}\left( x\right) \geq \frac{c}{V\left( x,\delta r_{1}\right) F_{1}}%
\sum_{i=c_{1}F_{1}}^{c_{2}F_{1}}\sum_{w\in D_{1}}c\mu \left( w\right) \geq c.
\end{equation*}%
\ Again from the definition of $u_{n}$ and $D_{1}\subset D,$ $%
c_{3}F_{1}<k<c_{4}F_{1},$ $x\in D_{1}$ we obtain 
\begin{equation}
\sum_{z\in D}\widetilde{p}_{k-m}^{B}\left( x,z\right) \mu \left( z\right)
\geq c  \label{d>c0}
\end{equation}%
or equivalently%
\begin{equation}
\sum_{z\in D}\widetilde{p}_{i}^{B}\left( x,z\right) \mu \left( z\right) \geq
c  \label{d>c}
\end{equation}%
if $x\in D_{1},\left( c_{3}-c_{2}\right) F_{1}<i<\left( c_{4}-c_{1}\right)
F_{1}.$

\item We will use the parabolic super mean value inequality in a new ball $%
B_{2}$ for $\widetilde{p}_{l}^{B}\left( x,y\right) $ with the same set of
constants $c_{i},$ hence with the same $\delta $ as well. Let $r_{2}=R/2$, $%
B_{2}=B\left( x,r_{2}\right) ,$ $D_{2}=B\left( x,\delta r_{2}\right) ,$ $%
F_{2}=F\left( x,r_{2}\right) .$ We apply $PSMV$ in $B_{2}$ and obtain that
for $c_{3}F_{2}<l<c_{4}F_{2},$ $y\in D_{2}$%
\begin{equation*}
\widetilde{p}_{l}^{B}\left( x,y\right) \geq \frac{c}{V\left( x,\delta
r_{2}\right) F\left( x,r_{2}\right) }\sum_{i=c_{1}F_{2}}^{c_{2}F_{2}}\sum_{z%
\in D_{2}}\widetilde{p}_{i}^{B}\left( x,z\right) \mu \left( z\right)
\end{equation*}%
if $\ $in addition $B_{2}\subset B\left( p,R\right) $. \ Let $x\in B\left( p,%
\frac{\delta R}{8}\right) $. \ This ensures that $B_{2}\subset B\left(
p,R\right) $ and $D_{2}\supset D$ and we obtain for $y\in B\left( x,\frac{%
\delta R}{4}\right) \subset D_{2},$ (and $B\left( x,\frac{\delta R}{4}%
\right) \subset B\left( p,\frac{\delta R}{2}\right) $ as well) that%
\begin{equation}
\widetilde{p}_{l}^{B}\left( x,y\right) \geq \frac{c}{V\left( x,\delta
r_{2}\right) F\left( x,r_{2}\right) }\sum_{c_{1}F_{2}}^{c_{2}F_{2}}\sum_{z%
\in D_{2}}\widetilde{p_{i}}^{B}\left( y,z\right) \mu \left( z\right) .
\label{pb>1}
\end{equation}%
In order to use $\left( \ref{d>c}\right) $ we require%
\begin{equation}
c_{2}F_{2}\geq \left( c_{4}-c_{1}\right) F_{1}  \label{f1}
\end{equation}%
and 
\begin{equation}
c_{1}F_{2}\leq \left( c_{3}-c_{2}\right) F_{1}.  \label{f2}
\end{equation}%
From the assumption $F\in W_{0}$ it follows that $\left( \ref{f1}\right) $
is satisfied if 
\begin{equation*}
c_{4}=c_{2}\left( 1+c_{F}4^{\beta ^{\prime }}\right)
\end{equation*}%
and $\left( \ref{f2}\right) $ is satisfied if%
\begin{equation*}
c_{1}=q\frac{\left( c_{3}-c_{2}\right) }{C_{F}}4^{-\beta }
\end{equation*}%
for any $0<q<1$. Finally%
\begin{equation*}
0<c_{4}-c_{1}=c_{2}\left( 1+c_{F}4^{\beta ^{\prime }}\right) -q\frac{\left(
c_{3}-c_{2}\right) }{C_{F}}4^{-\beta }<\varepsilon
\end{equation*}%
and $c_{1}<c_{2}<c_{3}<c_{4}$ can be ensured with the appropriate choice of $%
c_{2},c_{3}$ and $q$. Using $\left( \ref{f1}\right) $ and $\left( \ref{f2}%
\right) $ and $D_{2}\supset D$ the estimate in $\left( \ref{pb>1}\right) $
can be continued as follows:%
\begin{eqnarray*}
&&\widetilde{p}_{l}^{B}\left( x,y\right) \\
&\geq &\frac{c}{V\left( x,\delta r_{2}\right) F\left( x,r_{2}\right) }%
\sum_{c_{1}F_{2}}^{c_{2}F_{2}}\sum_{z\in D_{2}}\widetilde{p_{i}}^{B}\left(
y,z\right) \mu \left( z\right) . \\
&\geq &\frac{c}{V\left( x,\delta r_{2}\right) F\left( x,r_{2}\right) }%
\sum_{\left( c_{3}-c_{2}\right) F_{1}}^{\left( c_{4}-c_{1}\right)
F_{1}}\sum_{z\in D}\widetilde{p}_{i}^{B}\left( y,z\right) \mu \left(
z\right) .
\end{eqnarray*}%
Now we apply $\left( \ref{d>c}\right) $ to conclude to 
\begin{eqnarray*}
\widetilde{p}_{l}^{B\left( p,R\right) }\left( x,y\right) &\geq &\frac{c}{%
V\left( x,\delta r_{2}\right) F\left( x,r_{2}\right) }\sum_{\left(
c_{3}-c_{2}\right) F_{1}}^{\left( c_{4}-c_{1}\right) F_{1}}c \\
&\geq &\frac{c}{V\left( x,R\right) }\geq \frac{c}{V\left( x,f\left(
x,l\right) \right) },
\end{eqnarray*}%
where $c_{3}F_{2}\leq l\leq c_{4}F_{2}$ and $y\in B\left( p,\frac{\delta R}{4%
}\right) $. Finally let $S\geq 2R$%
\begin{equation}
\widetilde{p}_{l}^{B\left( p,S\right) }\left( x,y\right) \geq \widetilde{p}%
_{l}^{B\left( x,R\right) }\left( x,y\right) \geq \frac{c}{V\left( x,f\left(
x,l\right) \right) }  \label{PSL}
\end{equation}%
under the same conditions. \ Now choosing $\varepsilon ^{\prime }=\frac{c_{F}%
}{C_{F}}2^{-\beta }$ and $\delta ^{\prime }=\frac{\delta }{4}\left(
c_{3}c_{f}\right) ^{1/\beta ^{\prime }}$ $\left( \ref{PSL}\right) $ implies
that $PLE\left( F\right) $ (in the stronger form: $\left( \ref{sPLE}\right) $%
)%
\begin{equation}
\widetilde{p}_{l}^{B\left( p,S\right) }\left( x,y\right) \geq \frac{c}{%
V\left( x,f\left( x,l\right) \right) }
\end{equation}%
for $d\left( x,y\right) \leq \delta ^{\prime }f\left( p,l\right) ,l\leq
\varepsilon ^{\prime }F\left( p,S\right) .$
\end{enumerate}
\end{proof}

\subsection{Time comparison}

\label{sssum}In this subsection we summarize the results which lead to the
proof of $1\Longleftrightarrow 2\Longrightarrow 3\Longrightarrow 4$ in
Theorem \ref{tmainplussz} and we prove the return route from $%
4\Longrightarrow 2$ \ The equivalence of $1$ and $2$ is established by
Theorem \ref{tER1}, Corollary \ref{cg&H} and \ref{c*}, see also Remark \ref%
{rretutn}. The implication $2\Longrightarrow 3$ is given by Theorem \ref%
{tLDUE+} and \ref{tPLE} and $3\Longrightarrow 4$ is combination of Theorem %
\ref{tLDUE+} and \ref{tLNLE}, see also Remark \ref{rpr23}.

Now we prove $4\Longrightarrow 2$, the return route of Theorem \ref%
{tmainplussz}. Our task is to verify the implications in the diagram below
under the assumption $F\in W_{0}$ and $\left( p_{0}\right) ,\left( VD\right) 
$.%
\begin{gather}
\left. 
\begin{array}{c}
PMV_{1}\left( F\right) \\ 
PSMV\left( F\right)%
\end{array}%
\right\} \Longrightarrow \left. 
\begin{array}{l}
PMV_{\delta ^{\ast }}\left( F\right) \\ 
PSMV\left( F\right)%
\end{array}%
\right\} \Longrightarrow \left( H\right) \\
\left. 
\begin{array}{c}
PMV_{1}\left( F\right) \\ 
PSMV\left( F\right)%
\end{array}%
\right\} \Longrightarrow \left. 
\begin{array}{l}
DUE\left( F\right) \\ 
PLE\left( F\right) \\ 
\left( H\right)%
\end{array}%
\right\} \Longrightarrow \left. 
\begin{array}{c}
\rho v\simeq F \\ 
\left( H\right)%
\end{array}%
\right\} \Longrightarrow 
\begin{array}{c}
\left( TC\right) \\ 
\left( H\right)%
\end{array}%
\end{gather}

The heat kernel estimates are established as we indicated above. Now we deal
with proof of the elliptic Harnack inequality $\left( H\right) $ and the
time comparison principle $\left( TC\right) $.

\begin{theorem}
\label{tEHI}If $\Gamma $ satisfies $\left( p_{0}\right) ,\left( VD\right) $
and there is an $F\in W_{0}$ for which $PMV\left( F\right) $ and $PSMV\left(
F\right) $ are satisfied, then the elliptic Harnack inequality holds on $%
\left( H\right) $. \ 
\end{theorem}

We need an intermediate step, the parabolic mean value inequality for
smaller balls. We choose a particular set of constants $c_{i}$ subject some
restrictions coming from $PSMV$ and needed for later use.

\begin{lemma}
\label{cPMVs}\label{nodelta}If $(\Gamma ,\mu )$ satisfies $\left(
p_{0}\right) ,\left( VD\right) $ and $PMV_{1}\left( F\right) $ for an $F\in
W_{0}$, then for a given $\varepsilon ,\delta >0,0<$ $\delta ^{\ast }\leq 
\frac{1}{C_{F}}\varepsilon ^{\frac{1}{\beta ^{\prime }}}\frac{\delta }{2}$
there are $c_{1}<...<c_{4}$ such that $PMV_{\delta ^{\ast }}\left( F\right) $
holds for $\varepsilon $ and $c_{i}-s.$
\end{lemma}

\begin{proof}
We would like to derive $PMV_{\delta ^{\ast }}\left( F\right) $ for $c_{i}$
from $PMV_{1}\left( F\right) $ which holds for some other constants $a_{i}.$
We will apply $PMV_{1}\left( F\right) $ on the ball $B=B\left( x,\delta
R\right) $ and re-scale the time accordingly. We have $PMV_{\delta ^{\ast
}}\left( F\right) $ on $B\left( x,R\right) $ by%
\begin{eqnarray*}
\max_{\substack{ c_{3}F\left( x,R\right) \leq i\leq c_{4}F\left( x,R\right) 
\\ y\in B}}u_{i}\left( y\right) &\leq &\max_{\substack{ a_{3}F\left(
x,\delta R\right) \leq i\leq a_{4}F\left( x,\delta R\right)  \\ y\in B}}%
u_{i}\left( y\right) \\
&\leq &\frac{C}{\nu \left( D^{-}\right) }\sum_{j=a_{1}F\left( x,\delta
R\right) }^{a_{2}F\left( x,\delta R\right) }\sum_{y\in B}u_{j}\left(
z\right) \mu \left( z\right) \\
&\leq &\frac{C}{\nu \left( D^{-}\right) }\sum_{j=c_{1}F\left( x,R\right)
}^{c_{2}F\left( x,R\right) }\sum_{y\in B}u_{j}\left( z\right) \mu \left(
z\right)
\end{eqnarray*}%
if the inequalities $c_{1}<...<c_{4},a_{1}<...<a_{4},$%
\begin{eqnarray*}
a_{4}F\left( x,\delta R\right) &\geq &c_{4}F\left( x,R\right) \\
a_{3}F\left( x,\delta R\right) &\leq &c_{3}F\left( x,R\right) \\
a_{2}F\left( x,\delta R\right) &\leq &c_{2}F\left( x,R\right) \\
a_{1}F\left( x,\delta R\right) &\geq &c_{1}F\left( x,R\right)
\end{eqnarray*}%
are satisfied. \ We require in addition that $c_{4}\leq \varepsilon $ and $%
\frac{1}{C_{F}}\left( c_{3}-c_{2}\right) ^{\frac{1}{\beta ^{\prime }}}\frac{%
\delta }{2}\geq \delta ^{\ast }$. One can see that the following choice
satisfies these restrictions. Denote $p=C_{F}\left( \delta ^{\ast }\right)
^{\beta },q=c_{F}\left( \delta ^{\ast }\right) ^{\beta ^{\prime }}$. Let 
\begin{eqnarray*}
c_{4} &=&\varepsilon ,a_{4}=\frac{2q}{p}c_{4} \\
c_{3} &<&c_{4},a_{3}=qc_{3} \\
c_{2} &<&c_{3},a_{2}=\frac{1}{2}\min \left\{ pc_{2},a_{3}\right\} \\
c_{1} &=&\frac{1}{2}\min \left\{ \frac{a_{2}}{q},c_{2}\right\} ,a_{1}=qc_{1}.
\end{eqnarray*}%
Let us observe that $c_{1}$ can be arbitrarily small since $c_{4}\leq
\varepsilon $ and if the sub-solution is not given from an $m$ up to $%
a_{4}F\left( x,\delta R\right) $ it can be extended simply by $%
u_{i+m}=P_{i}^{B\left( x,R\right) }u_{m}$.
\end{proof}

\begin{proof}[Proof of Theorem \protect\ref{tEHI}]
Let us fix a set of constants $c_{1}<c_{2}<c_{3}<c_{4}=\varepsilon $ as in
Lemma \ref{cPMVs} and apply $PSMV\left( F\right) $ for them. Let us apply
Lemma \ref{cPMVs} for $\delta ^{\ast }$ to receive $PMV_{\delta ^{\ast
}}\left( F\right) $ on $B=B\left( x,R\right) .$ As a consequence for $%
D=B\left( x,\delta ^{\ast }R\right) ,u_{k}\left( y\right) =h\left( y\right) $
we obtain 
\begin{equation}
\max_{D}h\leq C\sum_{y\in D}h\left( y\right) .  \label{h1}
\end{equation}%
Similarly $PSMV\left( F\right) $ yields 
\begin{equation}
\min_{D}h\geq c\sum_{y\in D}h\left( y\right) .  \label{h2}
\end{equation}%
The combination of $\left( \ref{h1}\right) $ and $\left( \ref{h2}\right) $
gives the elliptic Harnack inequality \ for the shrinking parameter $\delta
^{\ast }.$ Finally $\left( H\right) $ can be shown using the standard
chaining argument along a finite chain of balls. The finiteness of the
number of balls follows from volume doubling via the bounded covering
principle.
\end{proof}

\begin{theorem}
\label{ttc}If $\left( p_{0}\right) $,$\left( VD\right) $ hold furthermore
there is an $F\in W_{0}$ for which $PMV\left( F\right) $ and $PSMV\left(
F\right) $ are satisfied, then $E\simeq F$ and $\left( TC\right) $ is true.
\end{theorem}

\begin{proposition}
\label{tERLE}Assume $\left( p_{0}\right) $ and $\left( VD\right) $ hold. If $%
PLE\left( F\right) $ holds for $F\in W_{0}$ hold, then then there is a $c>0$
such that for all $R>0,x\in \Gamma $%
\begin{equation*}
E\left( x,R\right) \geq cF\left( x,R\right) .
\end{equation*}
\end{proposition}

\begin{proof}
It follows from $PLE\left( F\right) $ that there are $c,C,1>\delta >\delta
^{\prime }>0,1>\varepsilon >\varepsilon ^{\prime }>0$ such that for all $%
x\in \Gamma ,R>1,$ $A=B(x,2R)$ and $n:\varepsilon ^{\prime }F\left(
x,R\right) <n<\varepsilon F\left( x,R\right) ,$ $r=\delta ^{\prime }R,$ $%
y\in B=B\left( x,r\right) $%
\begin{equation*}
\widetilde{P}_{n}^{A}(x,y)=P_{n}^{A}(x,y)+P_{n+1}^{A}(x,y)\geq \frac{c\mu (y)%
}{V(x,R)}.
\end{equation*}%
It follows for $F=\varepsilon F\left( x,R\right) ,F^{\prime }=\varepsilon
^{\prime }F\left( x,R\right) $ that 
\begin{align*}
E(x,2R)& =\sum_{k=0}^{\infty }\sum_{y\in B(x,2R)}P_{k}^{A}(x,y)\geq
\sum_{k=0}^{\infty }\sum_{y\in B}\frac{1}{2}\widetilde{P}_{k}^{A}(x,y) \\
& \geq \sum_{k=F^{\prime }}^{F}\sum_{y\in B}\frac{1}{2}\widetilde{P}%
_{k}^{A}(x,y)\geq c\frac{V(x,r)}{V(x,R)}F(x,R)\geq cF(x,R).
\end{align*}
\end{proof}

\begin{proposition}
\label{trrv<f}If $\left( p_{0}\right) ,\left( VD\right) $ hold and $%
DUE\left( F\right) $ holds for an $F\in W_{0}$, then there is a $C>0$ such
that for all $R>0,x\in \Gamma $%
\begin{equation*}
\rho \left( x,2R\right) v\left( x,2R\right) \leq CF\left( x,2R\right) .
\end{equation*}
\end{proposition}

The first step towards the upper estimate of $\rho v$ is to show an upper
estimate for $\lambda ^{-1}.$

\begin{proposition}
\label{pllf}If $(p_{0}),\left( VD\right) ,DUE\left( F\right) $ hold and $%
F\in W_{0}$, then there is a $c>0$ such that for all $R>0,x\in \Gamma $ 
\begin{equation}
\lambda (x,R)\geq cF^{-1}(x,R).  \label{llF}
\end{equation}
\end{proposition}

\begin{proof}
Assume that $C_{1}>1,2n=\left\lceil F(x,C_{1}R)\right\rceil ,$ $\ y,z\in
B=B(x,R)$. One can use \ 
\begin{equation*}
P_{2n}(y,z)=\sum_{w}P_{n}\left( y,w\right) P_{n}\left( w,z\right) \leq \sqrt{%
P_{2n}\left( y,y\right) P_{2n}\left( z,z\right) }
\end{equation*}%
and $DUE\left( F\right) $ to get%
\begin{equation*}
P_{2n}(y,z)\leq C\frac{\mu (z)}{\left( V(y,f(y,2n))V(z,f(z,2n))\right) ^{1/2}%
},
\end{equation*}%
(for the details see \cite{GT1}). From $\left( VD\right) $ and $F\in W_{0}$
it follows for $w=y$ or $z,$ $d(x,w)\leq R<C_{1}R=f(x,2n)$ 
\begin{equation*}
\frac{V(x,C_{1}R)}{V(w,C_{1}R)}\leq C,
\end{equation*}%
which results using $\left( p_{0}\right) $ that for all $n$ 
\begin{equation*}
P_{n}(y,z)\leq C\frac{\mu (z)}{V(x,f(x,n))}.
\end{equation*}%
If $\phi $ is the left eigenvector (measure) belonging to the smallest
eigenvalue $\lambda $ of $-\Delta ^{B}$ normalized to $(\phi 1)=1$, then 
\begin{align*}
(1-\lambda )^{2n}& =\phi P_{2n}^{B}1=\sum_{y,z\in B(x,R)}\phi
(z)P_{2n}^{B}(z,y)\leq \sum_{y\in B(x,R)}\frac{C\mu (y)}{\min_{z\in
B(x,R)}V(z,f(z,2n))} \\
& \leq C\max_{z\in B\left( x,R\right) }\left( \frac{R}{f\left( z,2n\right) }%
\right) ^{\alpha }=C\max_{z\in B\left( x,R\right) }\left( \frac{1}{C_{1}}%
\frac{f\left( x,2n\right) }{f\left( z,2n\right) }\right) ^{\alpha } \\
& \leq C\left( \frac{1}{C_{1}}C_{f}\right) ^{\alpha }\leq \frac{1}{2},
\end{align*}%
if $C_{1}=2C^{1/\alpha }C_{f}$. Using the inequality and $1-\xi \geq \frac{1%
}{2}\log \frac{1}{\xi }$ for $\xi \in \lbrack \frac{1}{2},1]$, where $\xi
=1-\lambda (x,R)$, one has 
\begin{equation*}
\lambda (x,R)\geq \frac{\log 2}{4n}\geq cF(x,C_{1}R)^{-1}>cF(x,R)^{-1}.
\end{equation*}
\end{proof}

\begin{proof}[Proof of Proposition \protect\ref{trrv<f}]
Let us recall from \cite{ter} that 
\begin{equation*}
\lambda \left( x,2R\right) \rho \left( x,R,2R\right) V\left( x,R\right) \leq
1
\end{equation*}%
in general, applying $\left( VD\right) $ and $\left( \ref{llF}\right) $
immediately yields the statement.
\end{proof}

\begin{proposition}
\label{trrv>F}If $\left( p_{0}\right) ,\left( VD\right) $ hold and $%
PLE\left( F\right) $ for an $F\in W_{0}$, then there is a $c>0$ such that
for all $R>0,x\in \Gamma $%
\begin{equation*}
\rho \left( x,R,2R\right) v\left( x,R,2R\right) \geq cF\left( x,2R\right)
\end{equation*}
\end{proposition}

\begin{proof}
The inequality $\left( \ref{E<rrv}\right) $ states that%
\begin{equation}
\rho \left( x,R,2R\right) v\left( x,R,2R\right) \geq \min_{z\in \partial
B\left( x,\frac{3}{2}R\right) }E\left( z,R/2\right) .
\end{equation}%
From Proposition \ref{tERLE} we know that%
\begin{equation*}
\min_{z\in \partial B\left( x,\frac{3}{2}R\right) }E\left( z,R/2\right) \geq
c\min_{z\in \partial B\left( x,\frac{3}{2}R\right) }F\left( z,R/2\right)
\end{equation*}%
and from $F\in W_{0}$ it follows that%
\begin{equation*}
\rho \left( x,R,2R\right) v\left( x,R,2R\right) \geq \min_{z\in \partial
B\left( x,\frac{3}{2}R\right) }F\left( z,R/2\right) \geq cF\left(
x,2R\right) .
\end{equation*}
\end{proof}

\begin{proof}[Proof of Theorem \protect\ref{ttc}]
From Proposition \ref{trrv<f} we have that $\rho v<CF$ which together with
Proposition \ref{trrv>F} yields that%
\begin{equation*}
\rho \left( x,R,2R\right) v\left( x,R,2R\right) \simeq F\left( x,2R\right) .
\end{equation*}%
Since $F\in W_{0}$ we have that $\rho v\in W_{0}$ and $\left( aD\rho
v\right) $ as well. From the conditions we have $\left( H\right) $ thanks to
Theorem \ref{tEHI} and by Theorem \ref{tER1} the Einstein relation follows:%
\begin{equation}
E\left( x,2R\right) \simeq \rho \left( x,R,2R\right) v\left( x,R,2R\right)
\simeq F\left( x,2R\right) .  \label{eccf}
\end{equation}%
Since $F\in W_{0}$ \ and $E\simeq F$ it follows that $E\in W_{0}$ which
includes $\left( TC\right) $ and of course $\left( wTC\right) $ as well and
the proof of $4\Longrightarrow 2$ of Theorem \ref{tmainplussz} is complete.
\end{proof}

\section{The parabolic Harnack inequality}

\setcounter{equation}{0}\label{spPH}In this section we will prove the
following extension of Theorem \ref{tph}.

\begin{theorem}
\label{tph+}If a weighted graph $(\Gamma ,\mu )$ satisfies $\left(
p_{0}\right) $, then the following statements are equivalent:

\begin{enumerate}
\item $\left( VD\right) $ hold \ and there is an $F\in W_{1}$ such that $%
g\left( F\right) $ is satisfied,

\item $\left( VD\right) ,\left( H\right) $ and $\left( \ast \right) $ holds
furthermore $E\in V_{1},$

\item $\left( VD\right) ,\left( H\right) $ and $\rho v\in V_{1},$

\item $\left( VD\right) $ and $UE\left( F\right) ,PLE\left( F\right) ,$ for
an $F\in W_{1}$ are satisfied,

\item $\left( VD\right) $ holds and there is an $F\in W_{1}$ such that $%
PMV\left( F\right) $ and $PSMV\left( F\right) $ are true,

\item there is an $F\in W_{1}$ such that the two-sided heat kernel estimate
hold: there are $C,\beta \geq \beta ^{\prime }>1,c>0$ such that for all $%
x,y\in \Gamma ,$ $n\geq d\left( x,y\right) $%
\begin{equation}
c\frac{\exp \left[ -C\left( \frac{F\left( x,d\right) }{n}\right) ^{\frac{1}{%
\beta ^{\prime }-1}}\right] }{V\left( x,f\left( x,n\right) \right) }\leq 
\widetilde{p}_{n}\left( x,y\right) \leq C\frac{\exp \left[ -c\left( \frac{%
F\left( x,d\right) }{n}\right) ^{\frac{1}{\beta -1}}\right] }{V\left(
x,f\left( x,n\right) \right) }  \label{tse2}
\end{equation}%
where $d=d\left( x,y\right) $,

\item there is an $F\in W_{1}$ such that $PH\left( F\right) $ holds.
\end{enumerate}
\end{theorem}

The equivalence of the statements $1-5$ are based on Theorem \ref%
{tmainplussz}. What is left is to incorporate $6$ and $7$. In this section
we show that the mean value inequalities for $F\in W_{1}$ are equivalent to
the parabolic Harnack inequality and to the two-sided heat kernel estimate $%
\left( \ref{tse2}\right) $. We will show the following implications:%
\begin{equation*}
\left. 
\begin{array}{c}
PMV \\ 
PSMV \\ 
\left( VD\right)%
\end{array}%
\right\} \Longrightarrow PH\Longrightarrow \left\{ 
\begin{array}{c}
DUE \\ 
DLE \\ 
PSMV%
\end{array}%
\right\} \Longrightarrow \left\{ 
\begin{array}{c}
\left( VD\right) \\ 
PMV \\ 
PSMV%
\end{array}%
\right. ,
\end{equation*}%
\begin{equation*}
\left. 
\begin{array}{c}
UE \\ 
PLE \\ 
\left( VD\right)%
\end{array}%
\right\} \Longleftrightarrow \left\{ 
\begin{array}{c}
UE \\ 
LE%
\end{array}%
\right.
\end{equation*}

\begin{theorem}
\label{tPM->PH}Assume $\left( p_{0}\right) $. Let $F\in W_{1}$, then the
following equivalence holds: 
\begin{equation*}
\left( VD\right) +PMV\left( F\right) +PSMV\left( F\right)
\Longleftrightarrow PH\left( F\right) .
\end{equation*}
\end{theorem}

\begin{remark}
We give direct proof of the statement instead of the ready alternative from 
\cite{De}. In Theorem 3.10 of \cite{De} the decomposition method shows that
for the parabolic Harnack inequality it is enough to show $UE\left( F\right) 
$ and $PLE\left( F\right) $ for the Dirichlet heat kernel on $B\left(
x,R\right) $. Since we know that $PMV\Longleftrightarrow UE$ and $%
PSMV\Longleftrightarrow PLE$ the proof is similar but works via the
Dirichlet heat kernel estimates. Here we prefer the direct route.
\end{remark}

\begin{proof}
The proof consist of several smaller steps.

1. First we show $PH\left( F\right) $ for Dirichlet solutions for a
particular set of constants. We choose $c_{1}<...<c_{4}$ and $\delta ^{\ast
}=\frac{1}{C_{F}}\left( c_{3}-c_{2}\right) ^{\frac{1}{\beta ^{\prime }}%
}\delta /2$ as in Lemma \ref{cPMVs}. \ Denote $\Phi ^{+}=\left[ c_{3}F,c_{4}F%
\right] \times B\left( x,\delta ^{\ast }R\right) $ and $\Phi ^{-}=\left[
c_{1}F,c_{2}F\right] \times B\left( x,\delta ^{\ast }R\right) $. Using Lemma %
\ref{cPMVs} we have for $\delta ^{\ast },$ $PMV_{\delta ^{\ast }}\left(
F\right) :$%
\begin{equation}
\max_{\Phi ^{+}}u\leq \frac{C}{\nu \left( \Phi ^{-}\right) }\sum_{\Phi
^{-}}u_{i}\left( z\right) \mu \left( z\right)  \label{pm1}
\end{equation}%
\ Let us choose $c_{6}>c_{5}>c_{4}.$ The parabolic super mean value
inequality $PSMV\left( F\right) $ with $\mathcal{D}^{+}=\left[ c_{5}F,c_{6}F%
\right] \times B\left( x,\delta ^{\ast }R\right) ,\mathcal{D}^{-}=\Phi ^{-}$
states that 
\begin{equation}
\min_{\mathcal{D}^{+}}\widetilde{u}\geq \frac{c}{\nu \left( \Phi ^{-}\right) 
}\sum_{\Phi ^{-}}\widetilde{u}_{i}\left( z\right) \mu \left( z\right) .
\label{pm2}
\end{equation}%
The combination of $\left( \ref{pm1}\right) $ and $\left( \ref{pm2}\right) $
results that 
\begin{equation}
\max_{D^{-}}u\leq C\min_{D^{+}}\widetilde{u}  \label{PHu}
\end{equation}%
which is the parabolic Harnack inequality for Dirichlet solutions for the
constants $c_{3}<c_{4}<c_{5}<c_{6},$ $\delta ^{\ast },$ in other words $%
D^{-}=\Phi ^{+},D^{+}=\mathcal{D}^{+}$.

2. Let us use the decomposition for an arbitrary solution $w\geq 0$ on $%
\mathcal{D}=\left[ 0,F\left( x,R\right) \right] \times B\left( x,R\right) $.
\ The nonnegative linear decomposition results in a Dirichlet solution $%
u\geq 0$ on $\mathcal{D}$ for which $u=w$ on $B\left( x,\delta ^{\ast
}R\right) $ and $u\leq w$ in general (for the details of the decomposition
method see \cite{De} proof of Theorem 3.10). \ Now we use $\left( \ref{PHu}%
\right) $%
\begin{equation*}
\max_{D^{-}}w=\max_{D^{-}}u\leq C\min_{D^{+}}u\leq C\min_{D^{+}}w.
\end{equation*}%
Which means that we have $PH\left( F\right) $ for all solutions and for the
given $c_{i}-$s and $\delta ^{\ast }$.

3. It is standard knowledge that if the (classical) parabolic Harnack
inequality holds for a set of constants $c_{i},\delta $, then it is true for
arbitrary set of constants as well (with an other $C$). This is the case if $%
F\in W_{1}$. The key is that $\beta ^{\prime }>1$ ensures that the time
dimension of the space-time cylinder shrinks faster than the space dimension
and the usual chaining argument can be applied.

4. The implication $PH\left( F\right) \Longrightarrow \left( VD\right) $ can
be seen along the lines of the classical proof (c.f. \cite{De}). \ First
from $PH\left( F\right) $ the diagonal upper and lower estimates are deduced
without change of the proof 
\begin{equation}
p_{m}\left( x,x\right) \leq \frac{C}{V\left( x,f\left( x,m\right) \right) }
\label{duef}
\end{equation}%
and%
\begin{equation}
\widetilde{p}_{n}\left( x,x\right) \geq \frac{c}{V\left( x,f\left(
x,n\right) \right) }.  \label{dlef}
\end{equation}%
The inequality for $n<cm$ 
\begin{equation}
p_{n}\left( x,x\right) \leq C\widetilde{p}_{m}\left( x,x\right)  \label{p<p}
\end{equation}%
can be obtained from $PH\left( F\right) $ with the proper choice of the
constants. Now let $n=\left\lfloor F\left( x,R\right) \right\rfloor
,m=\left\lceil F\left( x,A^{p}R\right) \right\rceil $, $p\geq 1$ and $A\geq
2 $ \ is chosen to satisfy $p>\frac{\beta }{\beta ^{\prime }}$ and $A>\left( 
\frac{C_{F}}{c_{F}}\right) ^{\frac{1}{p\beta ^{\prime }-\beta }}$. As a
result from $\left( \ref{duef}\right) ,\left( \ref{dlef}\right) $ and $%
\left( \ref{p<p}\right) $ one obtains $\left( VD\right) :$ 
\begin{equation*}
V\left( x,2R\right) \leq V\left( x,A^{p}R\right) \leq CV\left( x,R\right) .
\end{equation*}

5. The implication $PH\left( F\right) \Longrightarrow PSMV\left( F\right) $
is evident. As in step 4 we deduced $PH\left( F\right) \Longrightarrow
DUE\left( F\right) $ and $PMV\left( F\right) $ follows from Theorem \ref%
{tmainplussz}.
\end{proof}

\begin{remark}
The elliptic Harnack inequality is a direct consequence of the $F$-parabolic
one as it is true for the classical case.
\end{remark}

\begin{theorem}
\label{tLE}Assume that $\left( \Gamma ,\mu \right) $ satisfies $\left(
p_{0}\right) $ and $\left( VD\right) $. \ Then for any $F\in W_{1}$ 
\begin{equation*}
NDLE\left( F\right) \Longrightarrow LE\left( F\right) .
\end{equation*}
\end{theorem}

\begin{proof}
A \ modified version of Aronson's chaining argument gives the statement. The
proof uses varying radii for the chain of balls. We give the idea of the
modification (the other technical details can be seen following \cite{TD} or 
\cite{GT1}).

Denote $\delta $ the constant in $NDLE\left( F\right) $ and let $1>\delta
^{\prime }>0$ arbitrary. If $d\left( x,y\right) <\delta f\left( x,n\right) $
the statement follows from $NDLE,$ if $\delta ^{\prime }n\leq d\left(
x,y\right) \leq n$ it follows from $\left( p_{0}\right) $.

Assume that $\delta f\left( x,n\right) <d\left( x,y\right) <\delta ^{\prime
}n$. Consider a shortest path $\pi $ between $x$ and $y$, denote $d=d\left(
x,y\right) ,$ 
\begin{equation}
m=\left\lfloor \frac{n}{l\left( n,R,A\right) }\right\rfloor -1,  \label{mdef}
\end{equation}%
$R=f\left( x,n\right) $, $S=f\left( y,n\right) $, $A=B\left( x,d+R\right)
\cup B\left( y,d+S\right) $. \ Let $o_{1}=x$ and 
\begin{equation*}
r_{1}=\left\lceil \delta c_{0}f\left( o_{1},m\right) \right\rceil 
\end{equation*}%
and choose $o_{2}\in \pi :d\left( o_{1},o_{2}\right) =r_{1}-1$ and
recursively%
\begin{equation}
r_{i}=\left\lceil \delta c_{0}f\left( o_{i},m\right) \right\rceil 
\label{rdef}
\end{equation}%
and $o_{i}\in \pi :d\left( o_{i},o_{i+1}\right) =r_{i}-1$ and $d\left(
y,o_{i+1}\right) <d\left( y,o_{i}\right) .$ Denote $B_{i}=B\left(
o_{i},r_{i}\right) $. The iteration ends for the first \thinspace $j$ for
which $y\in B_{j}.$ From $F\in W_{0}$ and $z_{i+1}\in B_{i}$ it follows that 
\begin{equation}
c_{1}\leq \frac{f\left( z_{i},m\right) }{f\left( z_{i+1},m\right) }\leq C_{2}
\label{epere}
\end{equation}%
and the from triangle inequality it is evident that 
\begin{equation}
d\left( z_{i},z_{i+1}\right) \leq 2r_{i}+r_{i+1}\leq \left( 2+\frac{1}{c_{1}}%
\right) \delta c_{0}f\left( z_{i},m\right) .  \label{delt}
\end{equation}%
Here we specify $c_{0}=\left( 2+1/c_{1}\right) ^{-1}$. Let us recall the
definition of $l=l\left( n,d,A\right) $ 
\begin{equation}
\frac{n}{l}\geq \max_{z\in A}CE\left( z,\frac{d}{l}\right) ,  \label{ldefC}
\end{equation}%
taking the inverse one obtains:%
\begin{equation}
\min_{z\in A}f\left( z,\frac{1}{C}\frac{n}{l}\right) \geq \frac{d}{l}.
\end{equation}%
Let us choose $C$ in $\left( \ref{ldefC}\right) $ (using $F\in W_{1}$) such
that 
\begin{equation*}
f\left( o_{i},\frac{1}{C}\frac{n}{l}\right) \leq \delta c_{0}f\left( o_{i},%
\frac{n}{l}\right) =r_{i}.
\end{equation*}%
By the definition of $j$ 
\begin{equation*}
d>\sum_{i=1}^{j-1}r_{i}\geq \left( j-1\right) \frac{d}{l},
\end{equation*}%
consequently, $j-1\leq l.$%
\begin{equation*}
\left( \widetilde{P}_{m}\right) ^{j}\left( x,y\right) \geq \sum_{z_{1}\in
B_{0}}...\sum_{z_{j-1}\in B_{j-2}}\widetilde{P}_{m}\left( x,z_{1}\right) 
\widetilde{P}_{m}\left( z_{1},z_{2}\right) ...\widetilde{P}_{m}\left(
z_{j-1},y\right) .
\end{equation*}%
Now we use $NDLE$ to obtain%
\begin{eqnarray*}
\left( \widetilde{P}_{m}\right) ^{j}\left( x,y\right)  &\geq &\sum_{z_{1}\in
B_{0}}...\sum_{z_{j-1}\in B_{j-2}}\frac{c\mu \left( z_{1}\right) }{V\left(
x,f\left( x,m\right) \right) }...\frac{c\mu \left( y\right) }{V\left(
z_{j-1},f\left( z_{j-1},m\right) \right) } \\
&\geq &\min_{z_{2}\in B_{1}}...\min_{z_{j-1}\in B_{j-2}}c^{j-1}\frac{V\left(
o_{1},r_{1}\right) }{V\left( x,f\left( x,m\right) \right) }...\frac{V\left(
o_{j-1},r_{j-1}\right) }{V\left( z_{j-2},f\left( z_{j-2},m\right) \right) }%
\frac{\mu \left( y\right) }{V\left( z_{j-1},f\left( z_{j-1},m\right) \right) 
} \\
&\geq &\min_{z_{2}\in B_{1}}...\min_{z_{j-1}\in B_{j-2}}c^{j-1}\frac{\mu
\left( y\right) }{V\left( x,f\left( x,m\right) \right) }\frac{V\left(
o_{1},r_{1}\right) }{V\left( z_{2},f\left( z_{2},m\right) \right) }..\frac{%
V\left( o_{j-2},r_{j-2}\right) }{V\left( z_{j-1},f\left( z_{j-1},m\right)
\right) }.
\end{eqnarray*}%
If \ we use $\left( \ref{rdef}\right) ,\left( \ref{epere}\right) $ and $%
\left( VD\right) $ it follows that%
\begin{eqnarray}
\left( \widetilde{P}_{m}\right) ^{j}\left( x,y\right)  &\geq &\min_{z_{2}\in
B_{1}}...\min_{z_{j-1}\in B_{j-2}}\frac{c^{j-1}\mu \left( y\right) }{V\left(
x,f\left( x,m\right) \right) }\frac{V\left( o_{1},r_{1}\right) }{V\left(
z_{2},\frac{1}{\delta c_{0}c_{1}}r_{1}\right) }..\frac{V\left(
o_{j-2},r_{j-2}\right) }{V\left( z_{j-1},\frac{1}{\delta c_{0}c_{1}}%
r_{j-2}\right) }  \label{LEl} \\
&\geq &\frac{c^{j-1}\mu \left( y\right) }{V\left( x,f\left( x,m\right)
\right) }\left( c^{\prime }\right) ^{j-2}  \notag \\
&\geq &\frac{c\mu \left( y\right) }{V\left( x,f\left( x,n\right) \right) }%
\exp \left[ -C\left( j-1\right) \right]   \notag \\
&\geq &\frac{c\mu \left( y\right) }{V\left( x,f\left( x,n\right) \right) }%
\exp \left[ -Cl\right]   \label{cl}
\end{eqnarray}%
From Lemma 13.6 \ of \cite{GT1} we know that there is a $c>0$ such that 
\begin{equation*}
\widetilde{P}_{n}\geq c^{n-lm}\left( \widetilde{P}_{m}\right) ^{l}
\end{equation*}%
if $n\geq lm+l-1$. Let us note that from $\left( \ref{mdef}\right) $ it
follows that $n-lm+l\leq 3l$ which results in%
\begin{eqnarray*}
\widetilde{P}_{n}\left( x,y\right)  &\geq &c^{n-lm}\left( \widetilde{P}%
_{m}\right) ^{l}\left( x,y\right) \geq c^{\prime }\frac{c^{3l}\mu \left(
y\right) }{V\left( x,f\left( x,n\right) \right) }\exp \left( -Cl\right)  \\
&\geq &\frac{c\mu \left( y\right) }{V\left( x,f\left( x,n\right) \right) }%
\exp \left[ -C\left( \frac{F\left( x,d\left( x,y\right) \right) }{n}\right)
^{\frac{1}{\beta ^{\prime }-1}}\right] .
\end{eqnarray*}%
This finishes the proof of the lower estimate.
\end{proof}

\begin{proof}[Proof of Theorem \protect\ref{tph+}]
First of all we have seen that under conditions $\left( VD\right) +\left(
H\right) +\left( \ast \right) $ we have that%
\begin{equation*}
E\in W_{0}\text{ or }\rho v\in W_{0}
\end{equation*}%
and $\left( ER\right) $. If in addition $E\in V_{1},$ then both functions
belong to $W_{1}.$ On the other hand $\rho v\in V_{1}$ implies $\left(
aD\rho v\right) $ which is in the set of equivalent conditions $\left( \ast
\right) $, furthermore $W_{0}\cap V_{1}=W_{1}$ which shows the equivalence
of $2$ and $3$. Based on these observations the equivalence of $1-4$ and $5$
is established by Theorem \ref{tmainplussz}. The equivalence $5$ and $7$ is
given in Theorem \ref{tPM->PH}.

The implication $4\Longrightarrow 6$ follows from Theorem \ref{tmainplussz}
and \ref{tLE}. The reverse implication with respect to the upper estimate is
also covered by Theorem \ref{tmainplussz} as well. For any $F\in W_{1}$ $%
LE\left( F\right) $ implies $\left( VD\right) $. This can be seen exactly as
it is proved for $F\left( x,R\right) =R^{2}.$ \ The proof of $UE\left(
F\right) +LE\left( F\right) \Longrightarrow PLE\left( F\right) $ can be
reproduced following the steps of the proof of Lemma 8.3 in \cite{TD}. This
shows the equivalence of $4$ and $6$ and proves the whole statement.
\end{proof}

\begin{remark}
Let us note that we can prove a slightly better upper and lower estimates
(which are, in fact, equivalent to the ones presented). \ Denote $d=d\left(
x,y\right) $. \ Following the proof of the upper estimate in \cite{t2} (see
the proof of Theorem 3.14 and Remark 3.4) one can see that%
\begin{equation}
p_{n}\left( x,y\right) \leq \frac{C\exp \left[ -ck_{y}\left( n,\frac{1}{2}%
d\right) \right] }{V\left( x,f\left( x,n\right) \right) }+\frac{C\exp \left[
-ck_{x}\left( n,\frac{1}{2}d\right) \right] }{V\left( y,f\left( y,n\right)
\right) }.  \label{a}
\end{equation}%
The intermediate estimate $\left( \ref{cl}\right) $ gives a stronger lower
bound:%
\begin{equation}
\widetilde{p}_{n}\left( x,y\right) \geq \frac{c}{V\left( x,f\left(
x,n\right) \right) }\exp \left[ -Cl\left( x,n,A\right) \right] ,  \label{b}
\end{equation}%
where $A=B\left( x,d\left( x,y\right) +f\left( x,n\right) \right) \cup
B\left( y,d\left( x,y\right) +f\left( y,n\right) \right) ,n\geq d\left(
x,y\right) \,$.
\end{remark}

\begin{remark}
\label{rend}It is not immediate, but is elementary to deduce from $\left( %
\ref{a}\right) $ and $\left( \ref{b}\right) $ a spacial case of Theorem \ref%
{tph+} if%
\begin{equation*}
E\left( x,R\right) \simeq F\left( R\right)
\end{equation*}%
or%
\begin{equation*}
\rho \left( x,R,2R\right) v\left( x,R,2R\right) \simeq F\left( R\right) .
\end{equation*}%
Such a result is presented in \cite{t2}. The key observation is that under $%
\left( p_{0}\right) ,\left( VD\right) ,\left( H\right) $ the condition $%
\left( E\right) $ implies $\beta ^{\prime }>1$. The statements $1-5$ and $7$
of Theorem \ref{tph+} are immediate, the two-sided heat kernel estimate%
\begin{equation}
c\frac{\exp \left[ -Cm\left( n,d\left( x,y\right) \right) \right] }{V\left(
x,f\left( n\right) \right) }\leq \widetilde{p}_{n}\left( x,y\right) \leq C%
\frac{\exp \left[ -cm\left( n,d\left( x,y\right) \right) \right] }{V\left(
x,f\left( n\right) \right) }  \label{tsgem}
\end{equation}%
needs some preparation (here $f\left( n\right) $ is the inverse of $F\left(
R\right) $ again). It follows from $\left( \ref{a}\right) $ and $\left( \ref%
{b}\right) $ and from the fact that for any fixed $C_{i}>0,x\in \Gamma $%
\begin{equation*}
k_{x}\left( C_{1}n,C_{2}R\right) \simeq l_{x}\left( C_{3}n,C_{4}R\right)
\simeq m\left( C_{5}n,C_{6}R\right) .
\end{equation*}%
In the very particular case when $E\left( x,R\right) \simeq R^{\beta }$ one
recovers from $\left( \ref{tsgem}\right) $ the sub-Gaussian estimate $\left( %
\ref{sgebeta}\right) $ which is usual for the simplest fractal like graphs.
\end{remark}

\textbf{List of the main conditions}

\begin{equation*}
\begin{array}{lllll}
\text{shortcut} &  & \text{equation} &  & \text{name} \\ 
\left( p_{0}\right) &  & \left( \ref{p0}\right) &  & \text{controlled
weights condition} \\ 
\left( VD\right) &  & \left( \ref{vdef}\right) &  & \text{volume doubling
property} \\ 
\left( TC\right) &  & \left( \ref{TC}\right) &  & \text{time comparison
principle} \\ 
\left( ER\right) &  & \left( \ref{ER}\right) &  & \text{the Einstein relation%
} \\ 
\left( MV\right) &  & \left( \ref{MV}\right) &  & \text{mean value inequality%
} \\ 
DUE\left( E\right) &  & \left( \ref{LDUE}\right) &  & \text{diagonal upper
estimate} \\ 
DLE\left( F\right) &  & \left( \ref{dlef}\right) &  & \text{diagonal lower
estimate} \\ 
g\left( F\right) &  & \left( \ref{gUE}\right) +\left( \ref{gLE}\right) &  & 
\text{bounds on }g \\ 
\left( H\right) &  & \left( \ref{H}\right) &  & \text{elliptic Harnack
inequality} \\ 
UE\left( F\right) \text{ \ } &  & \left( \ref{UE}\right) \text{ \ \ } &  & 
\text{upper estimate w.r.t. }F\text{ \ \ } \\ 
PLE\left( F\right) &  & \left( \ref{PLE0}\right) &  & \text{particular lower
estimate} \\ 
NDLE\left( F\right) \text{ \ } &  & \left( \ref{NDLE}\right) \text{ \ \ } & 
& \text{near\ diagonal lower estimate} \\ 
LE\left( F\right) \text{ \ } &  & \left( \ref{LE}\right) \text{ \ \ } &  & 
\text{lower estimate} \\ 
PMV\left( F\right) \text{ \ } &  & \left( \ref{PMV}\right) \text{ \ \ } &  & 
\text{parabolic mean value inequality\ \ } \\ 
PSMV\left( F\right) &  & \left( \ref{minPH}\right) &  & \text{parabolic
super mean value inequality} \\ 
\left( E\right) ,\left( \rho v\right) &  & \left( \ref{E}\right) \left( \ref%
{rrvA}\right) \text{\ \ } &  & E\left( x,R\right) \text{\ or }\rho v\text{
is uniform in \ \ }x \\ 
\left( aD\rho v\right) \text{ \ } &  & \left( \ref{adrrv}\right) \text{ \ \ }
&  & \text{anti-doubling for }\rho v \\ 
\begin{array}{l}
PH\left( F\right) \\ 
RLEF%
\end{array}
&  & 
\begin{array}{l}
\left( \ref{PMV}\right) \\ 
\left( \ref{RLE}\right)%
\end{array}
&  & 
\begin{array}{l}
\text{parabolic Harnack inequality} \\ 
\text{resistance lower estimate}%
\end{array}%
\end{array}%
\end{equation*}

\end{document}